\documentclass[a4paper,11pt]{article}
\usepackage{amssymb}

\usepackage{graphicx}

\newcommand{\R}{\mathbb{R}}

\newcommand{\N}{\mathbb{N}}

\newcommand{\cuad}{{\sqcap\kern-.68em\sqcup}}

\newcommand{\norm}[1]{\|#1\|}

\newtheorem{theorem}{Theorem}[section]
\newtheorem{proposition}{Proposition}[section]

\newtheorem{lemma}{Lemma}[section]
\newtheorem{corollary}{Corollary}[section]
\newtheorem{remark}{Remark}[section]
\newcommand{\bremark}{\begin{remark} \em}
\newcommand{\eremark}{\end{remark} }

\usepackage{cite}
\begin{document}

\begin{center}{\bf  \Large   Fast decaying and slow decaying solutions of Lane-Emden \\[1.5mm]

equations  involving nonhomogeneous potential

 }\bigskip\medskip

 {\small  Huyuan Chen\footnote{chenhuyuan@yeah.net}
}
 \bigskip

{\small   Department of Mathematics, Jiangxi Normal University,\\
Nanchang, Jiangxi 330022, PR China\\[6mm]

 {\small   Xia Huang \footnote{xhuang1209@gmail.com  }\quad and\quad  Feng Zhou\footnote{ fzhou@math.ecnu.edu.cn}
} \bigskip

Center for PDEs, School of Mathematical Sciences, East China Normal University,\\
Shanghai Key Laboratory of PMMP,
Shanghai 200241, PR China \\[6mm]

}

\begin{abstract}
Our purpose of this paper is to study  positive solutions of Lane-Emden equation
 $$ -\Delta  u= V u^p\quad {\rm in}\quad \R^N\setminus\{0\} $$
disturbing by a non-homogeneous potential  $V$ when $p\in (\frac{N}{N-2}, p_c)$, where $p_c$ is the  Joseph-Ludgren exponent. We construct a sequence of fast decaying solutions and slow decaying solutions with appropriated restrictions for $V$.
\end{abstract}

  \end{center}

\setcounter{equation}{0}
\section{Introduction}

 Our concern in this paper is to consider fast decaying solutions of weighted Lane-Emden equation in punctured domain
\begin{equation}\label{eq 1.1}
 \left\{
\begin{array}{lll}
 \displaystyle  -\Delta  u= V u^p\quad    &{\rm in}\quad   \R^N\setminus\{0\},\\[2mm]
 \displaystyle  \quad\ \ u>0\quad    &{\rm in}\quad   \R^N\setminus\{0\},
 \end{array}\right.
\end{equation}
where $p>1$, $N\ge3$ and the potential $V$ is a locally H\"{o}lder continuous function in $\R^N\setminus\{0\}$.

When $V=-1$, the nonlinear term is known as an absorption and problem (\ref{eq 1.1})
admits a positive solution $c_p|x|^{-\frac{2}{p-1}}$ for any $p\in(0,\frac{N}{N-2})$ and
Brezis-Veron in \cite{BV} showed that it has no positive solution when $p\ge \frac{N}{N-2}$.
This type of nonexistence in the super critical case could also done in \cite{VV}.
While the isolated singularities of this elliptic problems in punctured domain subject to Dirichlet boundary condition are well studied in  \cite{GV,MV2,MV3,MV4,V0} and a survey in \cite{V}.

When $V=1$, equation (\ref{eq 1.1}) is well known as Lane-Emden-Fowler equation
\begin{equation}\label{ho}
  -\Delta  u=   u^p\quad    {\rm in}\quad   \R^N\setminus\{0\},
\end{equation}
which has been extensively studied in the last decades. When $p\le \frac{N}{N-2}$, problem (\ref{ho}) has no positive solution, see the reference \cite{BP} and for $p>\frac{N}{N-2}$, problem (\ref{ho})
always has a singular solution $w_p(x)=c_p|x|^{-\frac{2}{p-1}}$ with
\begin{equation}\label{cp}
 c_p=\left(\frac{2}{p-1}(N-2-\frac{2}{p-1})\right)^{\frac1{p-1}}.
\end{equation}
When $p\in(\frac{N}{N-2},\frac{N+2}{N-2})$,
positive isolated singular solutions of the problem (\ref{ho}) have the following structure:

\noindent (i) {\it a sequence $k$-fast decaying solutions $w_{k}$ with $k>0$ such that
$$
 \lim_{|x|\to0^+}w_{k}(x)|x|^{\frac{2}{p-1}}=c_p\quad{and}~\lim_{|x|\to+\infty}w_{k}(x)|x|^{N-2}=k.
$$

Here a solution is called $k$-fast decaying if $\displaystyle  \lim_{|x|\to+\infty}w_{k}(x)|x|^{N-2}=k$. }

\noindent(ii){\it a slow decaying solution $w_{p}(x)=c_p|x|^{-\frac{2}{p-1}}$ and $\displaystyle  w_{p}=\lim_{k\to+\infty}w_{k}$.} \\[1mm]
Furthermore, the fast decaying solution $w_k$ could be written by
$$w_{k}(x)=|x|^{-\frac{2 }{p-1}}\bar w_p(-\ln|x|+b_p ^{-1} (\ln k-\ln c)), $$
where $b_p = N-2-\frac{2}{p-1}>0$, $c>0$ is independent of $k$ and $\bar w_p(\cdot)$ is a positive and bounded function independently of $k$.
Assume that  $t=-\ln|x|+b_p ^{-1} (\ln k-\ln c)$, then the function $\bar w_p$ satisfies
\begin{equation}\label{eq 2.3-p}
  \left\{
\begin{array}{lll}
  \bar w_p'' -\left(N-2-\frac{4}{p-1}\right)  \bar w_p'-c_p^{p-1} \bar w_p+  \bar w_p^p=0\quad{\rm in}\quad \R,\\[2.5mm]
   \displaystyle
  \bar w_p(-\infty)=0\quad{\rm and}\quad  \bar w_p(+\infty)=c_p.
 \end{array}\right.
\end{equation}
To be convenient for the analyze, let us denote
\begin{equation}\label{pc}
 p_c=1+\frac4{N-4+2\sqrt{N-1}}\in\left(\frac{N}{N-2},\,\frac{N+2}{N-2}\right),
\end{equation}
which is the Joseph-Ludgren exponent, note that $\bar w_p$ is increasing for $p\in (\frac{N}{N-2},\ ,p_c]$ and for
$p\in (p_c, \frac{N+2}{N-2})$, $\bar w_p$ is oscillating as $t\to+\infty$, more information could be seen in Section 2.
For the supercritical case that $p\ge \frac{N+2}{N-2}$, problem (\ref{ho}) has been studied  \cite{DPM,DDG,MP1}. In particular, the authors
in \cite{DDG} obtained a sequence of fast decay solutions of (\ref{ho}) with $p>\frac{N+2}{N-2}$ in an exterior domain. 

During the last years there has been a renewed and increasing interest in the study of the semilinear elliptic equations
 with potentials, motivated by great applications in mathematical fields and physical fields,
 e.g. the well known scalar curvature equation in the study of Riemannian geometry,
 the scalar field equation for standing wave of nonlinear Schr\"{o}dinger and Klein-G\"{o}rden equations, the Matukuma equation, see a survey \cite{L2,N} and    more references on decaying solutions at infinity see \cite{CGS,CL,CPZ,DLY}.
For Lane-Emden equation (\ref{eq 1.1}) involving nonhomogeneous potential $V(x)=|x|^{\alpha_0}(1+|x|)^{\beta-\alpha_0},$
the authors in \cite{B,BP} showed the nonexistence provided $\beta>-2$ and $p\le \frac{N+\beta}{N-2}$, also see \cite[Theorem 3.1]{AS}.
In \cite{CFY}, the infinitely many positive solutions of problem (\ref{eq 1.1}) are constructed for
 $p\in(\frac{N+\beta}{N-2},\frac{N+\alpha_0}{N-2})\cap (0,\,+\infty)$ with $\alpha_0\in(-N,+\infty)$ and $\beta\in(-\infty,\alpha_0)$,
 by dealing with the distributional solutions of
\begin{equation}\label{eq 1.2}
 -\Delta u=Vu^p+\kappa\delta_0\quad {\rm in}\quad\R^N,
\end{equation}
where $k>0$, $\delta_0$ is a Dirac mass at the origin and $p=\frac{N+\alpha_0}{N-2}$ is the critical exponent named Serrin exponent,  the value for problem (\ref{eq 1.2}) with recoverable isolated singularities. Compared to the case $V\equiv1$, problem (\ref{eq 1.1}) would have totally different isolated singular solution structure for the super critical case $p\ge \frac{N+\alpha_0}{N-2}$, due to
the behavior of potential at infinity.

Our interest of this paper is to classify the fast decaying and slow decaying solutions of problem (\ref{eq 1.1})
for the supercritical case and involving general potential $V$.
Here, we say that {\it $u\in C^2(\R^N\setminus \{0\})$ is  a $\nu$-fast decaying solution
 if $u$ pointwisely verifies (\ref{eq 1.1}) in $\R^N\setminus\{0\}$ and has the asymptotic behavior at infinity
$$\lim_{|x|\to+\infty} u(x)|x|^{N-2}=\nu\quad{\rm for}\ \ \nu>0.$$}

Assume that the potential function $V$ is  H\"{o}lder continuous and satisfies the following conditions:
\begin{itemize}
\item[$(\mathcal{V}_0)$] 
$(i)$   near the origin,
\begin{equation}\label{v2-0}
 |V(x)-1|\le c_0|x|^{\tau_0}\quad{\rm for}\ \, x\in B_1(0),
\end{equation}
for some $c_0>0$ and $\tau_0>0$;\\[1mm]
$(ii)$ global control,
\begin{equation}\label{v1}
 0\le  V(x) \le c_\infty (1+|x|)^\beta\quad {\rm for}\quad |x|>0,
\end{equation}
where $c_\infty\geq 1$ and  $\beta\in\R$.

\end{itemize}

\medskip

Our main result is the following, which states the existence of fast decaying solutions of (\ref{eq 1.1}).

\begin{theorem}\label{teo 1}
Assume that $p_c$ is given by (\ref{pc}),  $p\in\left(\frac{N}{N-2},\,  p_c\right)$, the potential function $V$ verifies $(\mathcal{V}_0)$ with  $\tau_0$ and $\beta$ verifying
\begin{equation}\label{v00}
 \tau_0>\tau^*_p\quad{\rm and}\quad \beta< (N-2)p-N,
 \end{equation}
 where 
\begin{equation}\label{v001}
\tau^*_p=\Big(\frac{2}{p-1}-\frac{N-2}{2}\Big)-\sqrt{\Big(\frac{2}{p-1}-\frac{N-2}{2}\Big)^2-2\Big(N-2-\frac{2}{p-1}\Big)}.
 \end{equation}

Then there exists  $\nu_0>0$ such that   for any $\nu\in(0,\nu_0]$,   problem (\ref{eq 1.1})  has a  $\nu$-fast decaying  solution $u_\nu$,
which has singularity at the origin as
\begin{equation}\label{o1}
 \lim_{|x|\to0}u_\nu(x)|x|^{\frac{2}{p-1}}=c_{p},
\end{equation}
where $c_p$ is given in (\ref{cp}).

Furthermore, the mapping $\nu\in(0,\nu_0]\mapsto  u_\nu$ is  increasing, continuous and satisfies that
\begin{equation}\label{1.3}
 \lim_{\nu\to0} \norm{u_\nu}_{L^\infty_{loc}(\R^N\setminus\{0\})}=0.
\end{equation}
\end{theorem}

  We remark that  $\tau^*_p>0$ is well-defined 
due to   $ (\frac{2}{p-1}-\frac{N-2}{2} )^2>2 (N-2-\frac{2}{p-1} )$ for $p\in(\frac{N}{N-2},\, p_c)$.
Theorem \ref{teo 1} constructs a parameterized fast decaying solutions $u_\nu$  of (\ref{eq 1.1}) with $\nu\in I$ being an
   interval and more properties of the mapping $\nu\mapsto u_\nu$ are founded.
   Here the main difficulty is that the potential $V$ breaks the
scaling invariance of the equation. Moreover, due to the potential $V$, we can not restrict to search the symmetric solutions
by ODE's tools such as the phase analysis,  the variational method fails to apply due to the singularity at origin,.
First step of our method is to use the Schauder fixed point theorem to construct a solution $v_k$ of the problem
\begin{equation}\label{eq 00.1}
 -\Delta v= V(w_k+v)_+^p-w_k^p\quad{\rm in}\ \, \R^N\setminus \{0\},
\end{equation}
for $k>0$ sufficiently small and $w_k$ is the $k$-fast decaying solution of (\ref{ho}).
And then a $\tilde \nu_k$-fast decaying solution $\tilde u_{\nu_k}:=v_k+w_k$ of (\ref{eq 1.1})
is derived. However, the method of the Schauder fixed point theorem
fails to build the increasing mapping $k\mapsto \tilde\nu_k$. When $V$ is comparable to value 1, we note that
$v_k$ could be determined its sign, then motived by this observation, this solution could be
used as a barrier for the sequence  $v_n=\Gamma\ast(Vv_{n-1}^p)$ with initial data $v_0=w_k$,
and its limit is our desirable $\nu_k$-fast decaying solution of (\ref{eq 1.1}) and more properties of
the mapping $k\mapsto  \nu_k$ could be built, even for general $V$  by dividing it
as $V=(1+(V-1)_+)(1-(V-1)_-)$.

\smallskip

Our another interest of this paper is whether the parameter $\nu_0$ can be taken $+\infty$ in Theorem \ref{teo 1}, that is,
whether (\ref{eq 1.1}) has $\nu$-fast decaying solution  with $\nu\in(0,\,+\infty)$.  To this end,   we propose  the following assumptions
on the   potential $V$.

\begin{itemize}
\item[$(\mathcal{V}_1)$]
$(I)$ $V\geq 1$ in $\R^N\setminus\{0\}$ and there exist $\alpha_1\ge0$, $l_1>1$ such that
\begin{equation}\label{est 1}
 V(l_1^{-1}x)\geq l_1^{-\alpha_1} V(x), \quad\quad \forall\, x\in \R^N\setminus\{0\};
\end{equation}
$(II)$  $V\leq 1$ in $\R^N\setminus\{0\}$ ( i.e. $c_\infty=1$, $\beta=0$ in $(\mathcal{V}_0)$)  and there exist $\alpha_2\leq 0$, $l_2>1$ such that
 \begin{equation}\label{est 2}
V(l_2^{-1}x)\leq l_2^{-\alpha_2}V(x), \quad\quad \forall\, x\in \R^N\setminus\{0\}.
\end{equation}
\end{itemize}

\begin{theorem}\label{teo 2}
Assume that the potential function $V$ verifies $(\mathcal{V}_0)$ with   $\tau_0>0$ and $\beta$ verifying (\ref{v00})  and $p\in\left(\frac{N}{N-2},\,  p_c\right)$, where $p_c$ is given by (\ref{pc}).

 If $(\mathcal{V}_1)$ part $(I)$ or part $(II)$ holds,
 then for any $\nu\in(0,\,+\infty)$, problem (\ref{eq 1.1}) has a $\nu$-fast decaying solution $u_\nu$, which has singularity at the origin verifying (\ref{o1})
 and the mapping $\nu\in (0,\,\infty)\mapsto  u_\nu$ is increasing, continuous and (\ref{1.3})
 holds true.
 \end{theorem}

 Finally, our interest is to study the limit of $\{u_\nu\}_\nu$ as $\nu\to+\infty$ and we propose the following conditions on potential $V$.

 \begin{itemize}
\item[$(\mathcal{V}_\infty)$]
 Assume that $V$ is radially symmetric, decreasing with respect to $|x|$ and
\begin{equation}\label{v4}
\frac1\gamma |x| ^\alpha \le  V(x) \le  \gamma |x| ^\alpha\quad {\rm for}\quad |x|>1,
\end{equation}
where $\gamma>1$ and
\begin{equation}\label{poten 1}
(N-2)p_c-N-2<\alpha\leq 0.
\end{equation}
\end{itemize}

\begin{theorem}\label{teo 3}
Assume that $p_c$ is given by (\ref{pc}),  $p\in\left(\frac{N}{N-2},\,  p_c\right)$,
$V$ verifies $(\mathcal{V}_0)$ part $(i)$ with $\tau_0>0$,  $(\mathcal{V}_1)$ part $(II)$ and $(\mathcal{V}_\infty)$.
Let $u_\nu$ be a $\nu$-fast decaying solution of problem (\ref{eq 1.1}) with $\nu\in(0,+\infty)$ derived by Theorem \ref{teo 2}. Then the limit of $\{u_\nu\}_\nu$ as $\nu\to+\infty$ exists, denoting $\displaystyle  u_\infty=\lim_{\nu\to+\infty} u_\nu$, and $u_\infty$ is a solution of (\ref{eq 1.1}) verifying (\ref{o1}) and
\begin{equation}\label{1.3.2}
\frac1{c_1}\le  u_\infty(x)|x|^{\frac{2+\alpha}{p-1}}\le c_1,\quad |x|\ge 1,
\end{equation}
where $c_1>1$.
\end{theorem}

We see that the solution $u_\infty$ is no longer a fast decaying solution of (\ref{eq 1.1}) by the decay estimate (\ref{1.3.2}).
\emph{Here the solution $u_\infty$ of (\ref{eq 1.1}) is called as slow decay solution.}  Observe that the limit of $u_p=\lim_{k\to+\infty}w_k$, $k$-fast decay solution of (\ref{ho}),  behaviors
  as $c_p|x|^{-\frac{2}{p-1}}$ at infinity, and the fast decay solution $u_\nu$ behaviors as (\ref{1.3.2}), where
  $-\frac{2}{p-1}<-\frac{2+\beta}{p-1}.$
  This means $u_\infty>u_p$ for $|x|>r$ for some $r>0$, although  for any $\nu>0$, $u_\nu$ is derived by iterating the decreasing sequence
  $v_n=\mathbb{G}[Vv_{n-1}^p]$ with the initial data $v_0=w_k$ for some $k$ and $u_{\nu}\le w_k$. \smallskip

The rest of this paper is organized as follows. In section 2, we show qualitative properties of the solutions to
elliptic problem with homogeneous potential
and some basic estimates. Section 3 is devoted to build  fast decaying solutions of (\ref{eq 1.1}) by combining  Schauder fixed point theorem and iteration method.
Section 4 is devoted to consider the slow decaying solution as the limit of fast decaying solutions.

\setcounter{equation}{0}
\section{Preliminary}
\subsection{ODE analysis of Lane-Emden equation}

In this subsection, we recall some results of Lane-Emden equation by ODE analysis.  Denote
 a new independent variable $t=-\ln|x|$ and set $u(x)=|x|^{-\frac{2}{p-1}}\bar w_p(-\ln|x|)$, then the function $\bar w_p(t)$  verifies
\begin{equation}\label{eq 2.1}
  \bar w_p'' +a  \bar w_p'-c_p^{p-1} \bar w_p+  \bar w_p^p=0\quad{\rm in}\quad \R,
\end{equation}
where $a=\frac{4}{p-1}-N+2.$

Let $X(t)=\bar w(t)$ and $Y(t)=\bar w'(t)$, (\ref{eq 2.1}) can be rewritten as dynamic system
\begin{equation}\label{eq 2.1-ds}
 \left\{
\begin{array}{lll}
 \displaystyle  X'= Y\quad \    &{\rm in}\quad   \R ,\\[2mm]
 \displaystyle  Y'=-aY+c_p^{p-1}X-X^{p}    \ &{\rm in}\quad   \R.
 \end{array}\right.
\end{equation}
We see that system (\ref{eq 2.1-ds}) has two equilibrium points $(0,0)$ and $(c_p,0)$. Our aim
is to find a trajectory  $(X,Y)$ that  starts from point $(0,0)$ as $t\to-\infty$ and ends at
point $(c_p,0)$ as $t\to+\infty$.

In fact, the trajectory $(X,Y)$ is contained within the homocyclic orbit of the Hamiltonian system
$$ v'' -c_p^{p-1} v+  v^p=0\quad{\rm in}\quad \R,$$
then we conclude that $\sup\bar w_p\leq \sup v$.
The conservation of  Hamiltonian energy could  be given as following
$$E(t):=\frac12 v'(t)^2-\frac{c_p^{p-1} }2 v(t)^2+\frac1{p+1}v(t)^{p+1}=0.$$
Hence $v$ attains its supremum when $v'=0$, we get the upper bound $v^{p-1}<\frac{p+1}{2} c_p^{p-1}.$

On the other hand, the eigenvalues of linearizing system at $(0,0)$ to (\ref{eq 2.1-ds}) are
$$\lambda_1:=N-2-\frac{2}{p-1}>0>-\frac2{p-1}=:\lambda_2,$$
 for $t\to-\infty$, then we have that curve $(X,Y)$ goes out from point
$(0,0)$ along the direction $Y=\lambda_1X$ and then
$$\bar w_p(t)\sim ce^{\lambda_1 t}\quad {\rm as}\quad t\to-\infty.$$
Therefore, we have that
 \begin{equation}\label{2.1-1}
 \norm{\bar w_p}_{L^\infty(\R)}< (\frac{p+1}{2})^{\frac1{p-1}}\, c_p.
 \end{equation}
Observe that the eigenvalues at $(c_p,0)$ are given by zero points of
\begin{equation}\label{2.1}
H(\mu):= \mu^2+a \mu+(p-1)\, c_p^{p-1}.
\end{equation}
When $p\in(\frac{N}{N-2},p_c)$, where $p_c$ is given by (\ref{pc}), (\ref{2.1}) has two negative zero points of $H$:
$$\mu_1=\frac{1}{2}\left(-a+\sqrt{a^2-4(p-1)\, c_p^{p-1}}\,\right),\quad \mu_2=\frac{1}{2}\left(-a-\sqrt{a^2-4(p-1)\,c_p^{p-1}}\,\right),$$
then we have that
\begin{equation}\label{2.1--1}
 c_p-\bar w_p(t)\sim c_1e^{\mu_1 t}+c_2e^{\mu_2 t}\quad {\rm as}\quad t\to+\infty.
\end{equation}
For $p=p_c$, $H=0$ has two same roots $\mu_1=\mu_2=-\frac{a}{2}$, then
$c_p-\bar w_p(t)\sim (c_1+c_2t)e^{\mu_1 t}$ as $ t\to+\infty$.

Note that for $p\in(\frac{N}{N-2},\, p_c]$,
$\bar w_p$ is increasing and
\begin{equation}\label{2.5-2}
 \sup_{t\in\R}\bar w_p^{p-1}= \frac{2}{p-1}(N-2-\frac{2}{p-1})
\end{equation}
and then there exists $c>0$ such that
 \begin{equation}\label{e 2.2}
0< \bar w_p'\le c \bar w_p\quad{\rm in}\quad \R.
 \end{equation}
Hence, the fast decaying solutions $\{w_{k}\}_k$ of (\ref{ho}) have following properties:
$$
0< w_{k_1}< w_{k_2}< w_{p}\quad{\rm if}\quad 0<k_1< k_2<+\infty.
$$

When $p\in (p_c,\frac{N+2}{N-2})$, $H=0$ has two complex roots as
$$\mu_1=\frac{1}{2}\left(-a+{\rm i} \sqrt{-a^2+4(p-1)c_p^{p-1}}\right)\quad{\rm and}\quad \mu_2=\frac{1}{2}\left(-a-{\rm i}\sqrt{-a^2+4(p-1)c_p^{p-1}}\right),$$
where i is the unit imaginary number. Then $\bar w$ oscillates around $c_p$ and converges to $c_p$ with the rate
$$\limsup_{t\to+\infty} |\bar w_p(t)-c_p|e^{\frac a2 t}=c_0.$$

A crucial tool for the analysis is following.
\begin{proposition}\label{cr ex}
 $(i)$ For $p\in(\frac{N}{N-2},\, p_c]$, we have that
\begin{equation}\label{1.1}
 p\cdot \sup_{t\in\R}\bar w_p^{p-1}\le  \frac{(N-2)^2}{4},
\end{equation}
where $'='$ holds only for $p=p_c$.

 $(ii)$ For $p\in(p_c,\frac{N+2}{N-2})$, we have that
\begin{equation}
 p\cdot \sup_{t\in\R}\bar w_p^{p-1}> \frac{(N-2)^2}{4}.
\end{equation}
\end{proposition}
\noindent{\bf Proof.} One hand, when $p\in(\frac{N}{N-2}, \frac{N+2}{N-2})$, we have that
 $\frac{2}{p-1}\in (\frac{N-2}2, N-2)$,
 then
  $$\frac{2p }{p-1}(N-2-\frac{2}{p-1})\leq(\frac{N-2}{2})^2\quad{\rm holds ~if~ and~ only~ if}\quad  \frac{N}{N-2}<p\leq p_c.$$
On the other hand for $p\in (p_c,\frac{N+2}{N-2})$, we have that
$$\ \sup_{t\in\R}\bar w_p^{p-1}>\frac{2 }{p-1}(N-2-\frac{2}{p-1}),$$
combining (\ref{2.5-2}), the assertion holds and the proof is complete.\hfill$\Box$\medskip

\begin{lemma}\label{lm 2-1}
 Let $p\in (\frac{N}{N-2}, p_c]$ and $b_p=N-2-\frac{2}{p-1}$, then
\begin{equation}\label{2.2+1}
w_{k}(x)=|x|^{-\frac{2}{p-1}}\bar w_p(-\ln|x|+b_{p}^{-1} (\ln k-\ln d_0)),
\end{equation}
and for any $r\in(0,1]$, there exists $k_r=r^{b_{p}}$   such that  for $0<k\le k_r$,
\begin{equation}\label{2.2}
 w_{k}(x)\le c_1k r^{-\frac{2}{p-1}}(1+ |x|)^{2-N} \chi_{_{\R^N\setminus B_r(0)}}(x)+ c_p|x|^{-\frac{2}{p-1}}\chi_{_{B_r(0)}}(x)\quad{\rm for}\ \, \forall\, x\in\R^N\setminus\{0\},
\end{equation}
where $c_1>0$ is independent of $k,\, r$.
\end{lemma}
{\bf Proof.} By above phase plane analysis, we have (\ref{2.2+1})
just taking $t=-\ln|x|+b_{p}^{-1} (\ln k-\ln d_0)$.
From (\ref{2.1-1}), there exists $c>0$ such that
$\bar w_p(t)\leq c_1e^{b_pt}$ for any $t\leq 0$,
then if $-\ln|x|+b_{p}^{-1} (\ln k-\ln d_0)\leq0$, i.e.
$|x|\geq (k/d_0)^{b_{p}^{-1}}$, we have that
$$ w_{k}(x)\le c_1kr^{-\frac{2}{p-1}} (1+|x|)^{2-N}\quad \ \forall \,x\in \R^N\setminus B_r(0). $$
For $|x|<(k/d_0)^{b_{p}^{-1}}$, we have $\bar w_p\leq c_p$
and (\ref{2.2}) follows.\hfill$\Box$

 \begin{remark}
 Let $p\in (p_c,\frac{N+2}{N-2})$ and $b_p=N-2-\frac{2}{p-1}$, then
 for any $r\in(0,1]$, there exists $k_r=r^{b_{p}}$   such that  for $0<k\le k_r$,
$$
 w_{k}(x)\le c_1k r^{-\frac{2}{p-1}}(1+ |x|)^{2-N} \chi_{_{\R^N\setminus B_r(0)}}(x)+  \norm{\bar w_p}_{L^\infty} |x|^{-\frac{2}{p-1}}\chi_{_{B_r(0)}}(x),\quad\, \forall\, x\in\R^N\setminus\{0\},
$$
 where $\norm{\bar w_p}_{L^\infty}>c_p$.

 \end{remark}

\subsection{Basic estimate}

In this subsection, some estimates are introduced, which play important roles in
our construction of fast-decaying solutions for problem (\ref{eq 1.1}). Denote
$$\Gamma(x)=c_N|x|^{2-N},\quad \forall\, x\in\R^N\setminus\{0\},$$
which is the fundamental solution of $-\Delta \Gamma=\delta_0$ in $\R^N$
and $c_N>0$ is a normalized constant.

\begin{lemma}\label{lm 4.1}
 Let
 $$U_1(x)=|x|^{-2-\theta }\chi_{B_r(0)}(x),\quad U_2(x)=(1+|x|)^{-\tau}\quad{\rm and}\quad U_3(x)=|x|^{-\theta-2}(1+|x|)^{-\tau+\theta+2},$$
 where  $r\in(0,1/2)$ and $\tau>N>2+\theta$. Then there is $r^*>0$ small such that for $r\in(0, r^*]$,
 \begin{equation}\label{3.2}
  (\Gamma\ast U_1)(x)  \le  \frac{1}{\theta (N-2-\theta )} |x|^{-\theta }(1+|x|)^{2-N+\theta }\quad{\rm for} \ \, x\in\R^N
  \end{equation}
  and there exists $c>0$ such that
 \begin{equation}\label{3.3}
  (\Gamma\ast U_2)(x)   \le c  (1+|x|)^{2-N}\quad{\rm for} \ \, x\in\R^N
\end{equation}
and
\begin{equation}\label{3.4}
  (\Gamma\ast U_3)(x)   \le c  |x|^{-\theta}(1+|x|)^{2-N+\theta}\quad{\rm for} \ \, x\in\R^N.
\end{equation}
\end{lemma}
{\bf Proof.}  By direct computation, we have that
$$ (\Gamma\ast U_1)(x)=c_N \int_{B_{r}(0)}\frac{  |y|^{-2-\theta  }}{|x-y|^{N-2}}dy.$$
From the fact that $-\Delta (|x|^{-\theta})=\theta (N-2-\theta )|x|^{-2-\theta }$ in $\mathbb{R}^N\setminus\{0\}$, we can deduce
$$c_N\int_{\R^N}\frac{|y|^{-\theta -2} }{|x-y|^{N-2}}dy=\frac{1}{\theta (N-2-\theta )}|x|^{-\theta}.$$

For $x\in B_1(0)\setminus\{0\}$, we have that
 \begin{eqnarray*}
(\Gamma\ast U_1)(x)&\leq &  c_N\int_{\R^N}\frac{|y|^{-\theta -2} }{|x-y|^{N-2}}dz= \frac{1}{\theta (N-2-\theta )}|x|^{-\theta }.
\end{eqnarray*}
When $x\in \R^N\setminus B_1(0)$,
\begin{equation}\label{e 2-25}
  (\Gamma\ast U_1)(x) \leq c_N(|x|-r)^{2-N} \int_{B_r(0)} |y|^{-\theta -2}dy
\leq c_N 2^{N-2} |x|^{2-N} \int_{B_r(0)} |y|^{-\theta -2}dy,
\end{equation}
where
\begin{eqnarray*}
  \int_{B_r(0)} |y|^{-\theta -2}dy=|S^{N-1}|  r^{N-2-\theta}\to0\quad{\rm as}\ r\to0^+.
\end{eqnarray*}
Thus, (\ref{3.2}) holds thanks to $c_N 2^{N-2}|S^{N-1}|  r^{N-2-\theta}\le \frac{1}{\theta (N-2-\theta )}$ when  $r\leq r^*$. \smallskip

Next we show (\ref{3.3}). Note that
$$ (\Gamma\ast U_2)(x)=c_N \int_{\R^N}\frac{ (1+ |y|)^{-\tau }}{|x-y|^{N-2}}dy,$$
then $(\Gamma\ast U_2)(x)$ is bounded locally in $\R^N$ and so we only need to show the case $|x|\to+\infty$.

In fact, for $|x|>4$ large enough, there holds
\begin{eqnarray*}
\int_{\R^N}\frac{ (1+ |y|)^{-\tau }}{|x-y|^{N-2}}dy &\leq &\int_{B_{|x|/2}(0)}\frac{(1+|y|)^{-\tau} }{|x-y|^{N-2}}dy +\frac{1}{2}\int_{\R^N\setminus B_{|x|/2}(0)}\frac{|y|^{-\tau} }{|x-y|^{N-2}}dy\\
&\leq&\left(\frac{|x|}2\right)^{2-N}\int_{B_{|x|/2}(0)} (1+|y|)^{-\tau}  dy +|x|^{2-\tau}\int_{\R^N\setminus B_{1/2}(0)}\frac{|z|^{-\tau}}{|e_x-z|^{N-2}}dz\\
&\leq&c|x|^{2-N},
\end{eqnarray*}
where the last inequality holds thanks to the facts that $\tau>N$,
$$\int_{B_{|x|/2}(0)}(1+|y|)^{-\tau} dy\le \int_{\R^N}(1+|y|)^{-\tau} dy\quad {\rm and}\quad\int_{\R^N\setminus B_{\frac12}(0)} \frac{|z|^{-\tau} }{|e_x-z|^{N-2}}dz<+\infty.$$
Then (\ref{3.3}) holds true.\smallskip

Finally, we observe that
$$U_3(x)\le  2|x|^{-2-\theta}\chi_{B_1(0)}(x) + 2(1+|x|)^{-\tau},\quad x\in \R^N\setminus \{0\},$$
 then (\ref{3.4}) follows by (\ref{e 2-25}) with $r=1$ and (\ref{3.3}) directly.
\hfill$\Box$
\begin{corollary}\label{cr 2.2}
Assume that $\alpha\in (0,N)$,  $f$ is a nonnegative function satisfying that
$$ |f(x)|\le |x|^{-\theta}(1+|x|)^{\theta-\tau}\quad{\rm for}\quad |x|>0$$
with $\alpha<\theta<N$ and $\tau>N$. Then there exists $c>0$ such that
 \begin{equation}\label{2.4-0}
  \int_{\R^N}\frac{f(y)}{|x-y|^{N-\alpha}}\le c|x|^{-\theta+\alpha}(1+|x|)^{-N+\tau},\quad\forall\, x\in\R^N\setminus\{0\}.
 \end{equation}
\end{corollary}
{\bf Proof.} The same as the proof of Lemma \ref{lm 4.1}, we can obtain (\ref{2.4-0}).\hfill$\Box$

\begin{lemma}\label{lm 1.3}
Suppose that  $f\in L^1(\R^N)$ is a nonnegative function satisfying
$|f(x)|\le c|x|^{-\tau}$ for $|x|>r$,
with $\tau>N$ and some $r>0$, $c>0$. Then
 \begin{equation}\label{2.4}
 \lim_{x\to+\infty} (\Gamma\ast f)(x)|x|^{N-2}=c_N\int_{\R^N} f(x)dx.
 \end{equation}
\end{lemma}
{\bf Proof.} By the decay condition of $f$, we have that for any $\epsilon>0$, there exists $R>r_0$ such that for $R$ large,
$$\int_{B_R(0)} f(x) dx\geq (1-\epsilon) \norm{f}_{L^1(\R^N)}.$$
For $|x| \gg R$, there holds $(1-\epsilon) |x|^{2-N}\le |x-y|^{2-N}\le(1+\epsilon) |x|^{2-N}$ for $y\in B_R(0)$ and
\begin{eqnarray*}
(\Gamma\ast f)(x)   =  c_N \int_{B_R(0)} \frac{f(y)}{|x-y|^{N-2}}dy+ c_N\int_{\R^N\setminus B_R(0)} \frac{f(y)}{|x-y|^{N-2}}dy,
\end{eqnarray*}
which yields that for $|x|$ large,
\begin{eqnarray*}
 (1-\epsilon)  \norm{f}_{L^1(\R^N)}\le |x|^{N-2}\int_{B_R(0)} \frac{f(y)}{|x-y|^{N-2}}dy\le  (1+\epsilon)    \norm{f}_{L^1(\R^N)}
\end{eqnarray*}
and
\begin{eqnarray*}
\int_{\R^N\setminus B_R(0)} \frac{f(y)}{|x-y|^{N-2}}dy&\le& c \int_{\R^N\setminus B_R(0)} \frac{|y|^{-\tau}}{|x-y|^{N-2}}dy
\\&=&c  \int_{R\leq |y|< 2|x|} \frac{|y|^{-\tau}}{|x-y|^{N-2}}dy + c  \int_{|y|\geq 2|x|}\frac{|y|^{-\tau}}{|x-y|^{N-2}}dy\\
&\leq& c R^{2-N} \int _R^{2|x|} r^{N-1-\tau} dr + c\int^{+\infty}_{2|x|} r^{1-\tau} dr
\\&\leq & \frac{c}{N-\tau} R^{2-N} ((2|x|)^{N-\tau}-R^{N-\tau})-\frac{c}{2-\tau} (2|x|)^{2-\tau}.
\end{eqnarray*}
Passing to the limit as $\epsilon \to 0$ and letting $R\to+\infty$, we see that  $|x| \to +\infty$ and then (\ref{2.4}) holds.\hfill$\Box$

\begin{corollary}\label{cr 2.1}
Let $w_k$ be $k$-fast decaying solution of (\ref{ho}), then
$$c_N\int_{\R^N} w_k(x)^p dx=k.$$

\end{corollary}
{\bf Proof.} Since $\displaystyle \lim_{|x|\to0^+}w_k(x)|x|^{\frac{2}{p-1}}=c_p$ and
$\displaystyle  \lim_{|x|\to+\infty} w_k(x)|x|^{(N-2)}=k$,
 then Lemma \ref{lm 1.3} implies that for $p(N-2)>N$, which implies that
 $$k=\lim_{|x|\to+\infty} w_k (x)|x|^{ N-2}= c_N\int_{\R^N} w_k(x)^p dx,$$
which ends the proof.\hfill$\Box$

\begin{lemma}\label{lm 2.1}
Assume that  $a>0,\, b\in\R$,
then for $p\in (1,2]$,
$$(a+b)_+^{p} \le a^p+pa^{p-1} |b|+|b|^p;$$
for   $p>2$,
$$(a+b)_+^p\le a^p+pa^{p-1} |b|+2^{p} p(p-1) a^{p-2} b^2 +2^p|b|^p.$$
\end{lemma}
These are basic inequalities, here we omit the proof. 
\medskip

Finally, we introduce a comparison principle for general Hardy operator.

 \begin{lemma}\label{lm cp}
 Let $\Omega$ be a bounded $C^2$ domain containing the origin,
 $W$ be H\"{o}lder continuous locally in $\bar \Omega\setminus\{0\}$ such that $\displaystyle \lim_{|x|\to0} W(x)|x|^2=\mu $ with $\mu\in \Big(0, \frac{(N-2)^2}{4}\Big)$ and $W(x)\le  \frac{(N-2)^2}{4}|x|^{-2}$ in $\Omega$. Then the operator
\begin{equation}\label{lu}
 \mathcal{L}_W w:= -\Delta  w-W w
\end{equation}
verifies the following comparison principle in $\Omega$:

Assume that $f_1$, $f_2$ are two functions in $C^\gamma (\Omega\setminus\{0\})$ with $\gamma\in(0,1)$, $g_1$, $g_2$ are two continuous functions on $\partial \Omega$,
$$ f_1\ge f_2\quad {\rm in}\quad \Omega\setminus\{0\} \quad{\rm and}\quad  g_1\ge g_2\quad {\rm on}\quad \partial \Omega.$$
Let $u_i$ ($i=1,2$) be the classical solutions of
$$
 \arraycolsep=1pt\left\{
\begin{array}{lll}
 \displaystyle \mathcal{L}_W u = f_i\quad\
   &{\rm in}\ \,  {\Omega}\setminus \{0\},\\[1.5mm]
 \phantom{ L_\mu     }
 \displaystyle  u= g_i\quad  &{\rm   on}\ \, \partial{\Omega}.
 \end{array}\right.
$$
If $\displaystyle \liminf_{x\to0}u_1(x)|x|^{-\tau_-(\mu)} \ge \limsup_{x\to0}u_2(x)|x|^{-\tau_-(\mu)}$ holds, then $u_1\ge u_2 $ in $\Omega\setminus\{0\}.$

\end{lemma}
{\bf Proof.} From \cite[Theorem 1.1]{C}, we note that operator $\mathcal{L}_W$ has a positive solution $\Phi_W$ such that
$$\lim_{|x|\to0}\frac{\Phi_W(x)}{|x|^{\tau_-(\mu)}}=1.$$
Let $w=u_2-u_1$ be a solution of
 $$\arraycolsep=1pt\left\{
\begin{array}{lll}
 \displaystyle \mathcal{L}_W u \le 0\qquad
   {\rm in}\quad  {\Omega}\setminus \{0\},\\[1.5mm]
 \phantom{ L_\mu     }
 \displaystyle  u\le 0 \qquad  {\rm   on}\quad \partial{\Omega},\\[1.5mm]
 \phantom{   }
  \displaystyle \limsup_{x\to0}u(x)|x|^{-\tau_-(\mu)}\le0,
 \end{array}\right.$$
  then for any $\epsilon>0$, there exists $r_\epsilon>0$ converging to zero as $\epsilon\to0$ such that
 $w\le \epsilon \Phi_W$ on $\partial B_{r_\epsilon}(0)$.
Observe that
$w\le 0<\epsilon \Phi_W$ on  $\partial\Omega$,
then from \cite[Lemma 2.1]{CQZ}, we have that
$w\le \epsilon \Phi_W$ in $\Omega\setminus\{0\}$.
By the arbitrary of $\epsilon>0$, we have that $w\le 0$ in $\Omega\setminus\{0\}$.
\hfill$\Box$

\setcounter{equation}{0}
\section{Fast decaying solutions }

\subsection{Existence by Fixed point theory}
In this subsection, we give the proof of Theorem \ref{teo 1}.
We are looking for   a $k$-fast decaying solution $u$ of (\ref{eq 1.1}), with the division  
$$u=w_k+v,$$
where $w_k$ is the $k$ fast decaying solution of (\ref{ho}) and $v$ verifies that
\begin{equation}\label{eq 3.1}
 -\Delta v= V(w_k+v)_+^p-w_k^p\quad{\rm in}\ \, \R^N\setminus \{0\}.
\end{equation}
We will employ the Schauder fixed point theorem to obtain a solution of (\ref{eq 3.1}).  To this end, let us clarify 
the key value $\tau_p^*$.  In fact, the essential point in our following construction of fast decaying solutions 
is to find a $\theta_0\in [\frac{N-2}{2}, N-2)$ such that 
\begin{equation}\label{eq 3.1-001}
\theta_0(N-2-\theta_0)>pc_p^{p-1}= \left(2+ \frac{2}{p-1}\right)\left(N-2-\frac{2}{p-1}\right),
\end{equation}
which is possible since $\left(2+ \frac{2}{p-1}\right)\left(N-2-\frac{2}{p-1}\right)<\frac{(N-2)^2}{4}$ when $p\in(\frac{N}{N-2}, p_c)$,
and $\frac{(N-2)^2}{4}$ is maximum of $\theta_0(N-2-\theta_0)$ which is achieved at $\theta_0=\frac{N-2}{2}$.
 Now the point is to find the smallest $\tau>0$ such that 
\begin{equation}\label{3.1}
\Big(\frac{2}{p-1}-\tau\Big)\Big(N-2+\tau-\frac{2}{p-1}\Big)= \left(2+ \frac{2}{p-1}\right)\left(N-2-\frac{2}{p-1}\right).
\end{equation}
Direct computation shows that  for $p\in(\frac{N}{N-2}, p_c)$,  $\tau_p^*>0$ defined in (\ref{v001}) is the smallest zero of (\ref{eq 3.1-001}) and letting $\tau_p^\#= (\frac{2}{p-1}-\frac{N-2}{2} )+\sqrt{ (\frac{2}{p-1}-\frac{N-2}{2} )^2-2 (N-2-\frac{2}{p-1} )}$,  for any $\tau\in\left(\tau_p^*,  \tau_p^\#\right)$, we have that 
$$\Big(\frac{2}{p-1}-\tau\Big)\Big(N-2+\tau-\frac{2}{p-1}\Big)> \left(2+ \frac{2}{p-1}\right)\left(N-2-\frac{2}{p-1}\right).$$

 Now let us fix 
\begin{equation}\label{eq 3.1-002}
\tau_{1}=  \tau_p^*+\frac12\min\Big\{\tau_0-\tau_p^*,\, \frac{2}{p-1}-\frac{N-2}2\Big\} \ \, {\rm and}\ \, \theta_0=\frac{2}{p-1}-\tau_1,
\end{equation}
then $\tau_1 \in \left(\tau_p^*,  \tau_p^\#\right)\cap (0,\tau_0)$ and $\theta_0\in   \left[\frac{N-2}{2}, \frac{2}{p-1}\right)$ verifying 
\begin{equation}\label{eq 3.1-003}
\theta_0(N-2-\theta_0)>\frac{2p}{p-1}(N-2-\frac{2}{p-1}).
\end{equation}
 Finally, we denote $\tau_{2}=\tau_0-\tau_{1}>0.$

\begin{proposition}\label{pr 3.1}
Assume that $p_c$ is given by (\ref{pc}),  $p\in\left(\frac{N}{N-2},\,  p_c\right)$, the potential function $V$ verifies $(\mathcal{V}_0)$ with $\tau_0,\beta$ verifying
(\ref{v00}). Then there exist $k^*>0$ and $c>0$ such that for any $k\in(0,\, k^*)$, problem (\ref{eq 3.1}) has a classical solution $v_k$ such that
 \begin{equation}\label{4.2}
  |v_k(x)|\le ck |x|^{-\theta_0}(1+|x|)^{2-N+\theta_0},\quad\forall\, x\in\R^N\setminus\{0\}.
 \end{equation}
where    $\theta_0$ is defined in (\ref{3.1}).

\end{proposition}

\noindent{\bf Proof.} {\it Step 1: to show basic setting for applying the Schauder fixed point Theorem.}
Note that  for $p\in (\frac{N }{N-2}, p_c)$, there holds
 \begin{equation}\label{2.3}
\theta_0(N-2-\theta_0)>\frac{2p}{p-1}(N-2-\frac{2}{p-1}).
 \end{equation}

 Let $q_0\in \left(\frac{N}{N-1},\, \frac{N}{\theta_0+1}\right)$, we denote
\begin{equation}\label{dom 1}
\mathcal{D}_{\epsilon}:=\left\{v\in L^{q_0}(\R^N):\, |v(x)|\le \epsilon|x|^{-\theta_0}(1+|x|)^{2-N+\theta_0},\ \forall\, x\in \R^N\setminus\{0\} \right\},
\end{equation}
and
\begin{equation}\label{op l}
\mathcal{T}  v:=\Gamma\ast \left(V(w_k+v)_+^p-w_k^p \right),\quad\forall\, v\in \mathcal{D}_{\epsilon},
\end{equation}
where $\Gamma$ is the fundamental solution of $-\Delta$ in $\R^N$.\medskip

{\it Step 2:   to prove $\mathcal{T} \mathcal{D}_{\epsilon}\subset \mathcal{D}_{\epsilon}$ for $\epsilon, ~k>0$ small suitably.}  
\smallskip

{\it Case i: $p\in\left(\frac{N}{N-2}, p_c\right)\cap (1,2]$ (this happens when $N\geq 5$).} For given $\epsilon>0$ small and $v\in \mathcal{D}_{\epsilon}$, we have that
\begin{equation}\label{4.1}
 |\Gamma\ast \left(V(w_k+v)_+^p-w_k^p \right)|\le \Gamma\ast\left(|V-1|w_k^p+pVw_k^{p-1}|v| +V|v|^p  \right).
\end{equation}
As $\frac{p c_p^{p-1} }{\theta_0(N-2-\theta_0)}<1$ by (\ref{eq 3.1-003}), then there exists  $r^*\in(0,1]$  such that
 $ \frac{p c_p^{p-1} }{\theta_0(N-2-\theta_0)}\max_{B_{r^*}(0)} V<1$.
and we denote
$$\rho_0=1-\frac{p c_p^{p-1} }{\theta_0(N-2-\theta_0)}\max_{B_{r^*}(0)} V>0. $$

By (\ref{2.2}) and Lemma \ref{lm 4.1}, direct computation shows that
\begin{eqnarray*}
&&\Gamma\ast (|V-1|w_k^p)\\&\le& c_0c_p^pr^{\tau_2}\int_{B_r(0)}\frac{|y|^{\tau_1-\frac{2}{p-1}-2}}{|x-y|^{N-2}}dy +c_1^pr^{-\frac{2p}{p-1}}  c_\infty k^p \int_{\R^N}\frac{(1+|y|)^{p(2-N)+\beta}}{|x-y|^{N-2}}dy
\\&\le & c_2r^{ \tau_2} |x|^{-\theta_0}(1+|x|)^{2-N+\theta_0}+c_3 r^{-\frac{2p}{p-1}}   k^p  (1+|x|)^{2-N},
\end{eqnarray*}
where $c_2,\, c_3>0$ are independent of $k,\, r$.
Now we fix $r\in(0,\, r^*]$ such that
$ c_2r^{\tau_2}\le \frac18{\rho_0} $,
then for that $r$, there exists $k^*_1>0$ such that for $k\in(0,k^*_1)$,
\begin{equation}\label{choose k-2}
  c_3 k^pr^{-\frac{2p}{p-1}}\le  \frac{\rho_0} 8\epsilon.
\end{equation}
Therefore, we have that
$\Gamma\ast (|V-1|w_k^p)(x)\leq \frac{\rho_0}4\epsilon |x|^{-\theta_0}(1+|x|)^{2-N+\theta_0}$.
Moreover, we observe that
\begin{eqnarray*}
&&p\Gamma\ast (V w_k^{p-1}|v|)\\&\le& \epsilon p\left[\max_{B_{r^*}(0)} V\, c_p^{p-1} \int_{B_{r^*}(0)}\frac{ |y|^{-2-\theta_0}}{|x-y|^{N-2}}dy + c_{\infty}c_1^{p-1} k^{p-1}r^{-2} \int_{\R^N}\frac{|y|^{-\theta_0}(1+|y|)^{p(2-N)+\beta+\theta_0}}{|x-y|^{N-2}}dy\right]
\\&\le &\epsilon \left(\max_{B_{r^*}(0)} V \frac{pc_p^{p-1} }{\theta_0(N-2-\theta_0)}+c_{\infty}c_1^{p-1}p k^{p-1}r^{-2} \right)|x|^{-\theta_0}(1+|x|)^{2-N+\theta_0}
\\&\le &  (1-\frac{\rho_0}2)    \epsilon |x|^{-\theta_0}(1+|x|)^{2-N+\theta_0},
\end{eqnarray*}
where  $c_\infty,c_1>0$ are independent of $r,\, k$. Then for fixed $r$,  there exists $k^*_1>0$ such that for $k\in(0,k^*_1)$,
\begin{equation}\label{choose k-1}
c_{\infty}c_1^{p-1}p k^{p-1}r^{-2} < \frac{\delta_0}2.
\end{equation}

Furthermore, we note that
\begin{eqnarray*}
\Gamma\ast (V|v|^p )&\le&  c_\infty \epsilon^p \int_{\R^N}  \frac{|y|^{-\theta_0p}  (1+|y|)^{p(2-N)+\beta+\theta_0p}}{|x-y|^{N-2}}dy
\\&\le&  c_4\epsilon^p |x|^{-\theta_0}(1+|x|)^{2-N+\theta_0}
\le   \frac{\rho_0}4   \epsilon |x|^{-\theta_0}(1+|x|)^{2-N+\theta_0},
\end{eqnarray*}
where $c_4>0$ is independent of $k$, we choose
\begin{equation}\label{choose k-3}
 c_4\epsilon^{p-1}  \le \frac{\rho_0}4.
\end{equation}

As a consequence, for $\epsilon,r,k$ verifying (\ref{choose k-1})-(\ref{choose k-3}), we have that
$$ |\Gamma\ast \left(V(w_k+v)^p-w_k^p \right)|\le \epsilon |x|^{-\theta_0} (1+|x|)^{2-N+\theta_0},$$
that is to say,  $\mathcal{T} \mathcal{D}_{\epsilon}\subset \mathcal{D}_{\epsilon}$.\medskip

{\it Case ii: $p\in(\frac{N}{N-2}, p_c)\cap (2,+\infty)$.} In this case, (\ref{4.1}) should be replaced by
\begin{equation}\label{4.1-1}
 |\Gamma\ast \left(V(w_k+v)_+^p-w_k^p \right)|\le \Gamma\ast\left(|V-1|w_k^p+pVw_k^{p-1}|v|+2^pp(p-1)w_k^{p-2}v^2 +V|v|^p  \right).
\end{equation}
Here we only need to do the estimate for $\Gamma\ast(V w_k^{p-2}v^2)$ in addition. Indeed, as $w_k(x)\le c_p |x|^{-\frac{2}{p-1}}(1+|x|)^{\frac{2}{p-1}+2-N}$ for $k\le k_0$, we have that
\begin{eqnarray*}
\Gamma\ast (V w_k^{p-2}v^2 )&\le&  c_\infty c_p^{p-2} \epsilon^2\int_{\R^N}  \frac{|y|^{-\frac{2(p-2)}{p-1}-2\theta_0}  (1+|y|)^{\frac{2(p-2)}{p-1}+p(2-N)+\beta+2\theta_0}}{|x-y|^{N-2}}dy
\\&\le& c_5 \epsilon^2   |x|^{-\theta_0}(1+|x|)^{2-N+\theta_0}
\le   \frac{\rho_0}4   \epsilon |x|^{-\theta_0}(1+|x|)^{2-N+\theta_0},
\end{eqnarray*}
where $c_5>0$ is independent of $\epsilon$, $-\frac{2(p-2)}{p-1}-2\theta_0>-\theta_0-2$,   and
the last inequality holds if we choose
\begin{equation}\label{choose k-4}
  c_5\epsilon \le \frac{\rho_0}4.
\end{equation}
By our choice of  $\epsilon$, $r$ and $k$, we have that
$\mathcal{L} \mathcal{D}_{\epsilon}\subset \mathcal{D}_{\epsilon}$ for $p\in(\frac{N}{N-2}, p_c)\cap (2,+\infty)$.

\medskip
{\it Step 3: Applying Schauder fixed point theorem.}
Note that
 for $x\in \R^N\setminus\{0\}$,
$$h(x):= \left|  V(w_k(x)+v(x))_+^p-w_k^p(x) \right|\le c |x|^{-\theta_0-2}(1+|x|)^{p(2-N)+\beta},$$
and then by {\it Step 2} and  Corollary \ref{cr 2.2}, we have that
\begin{eqnarray*}
| \mathcal{T} v(x)|  \le c_N \int_{\R^N}\frac{h(y)}{|x-y|^{N-2}}dy
\le \epsilon |x|^{-\theta_0}(1+|x|)^{2-N+\theta_0}
\end{eqnarray*}
 and
\begin{eqnarray*}
|\nabla \mathcal{T} v(x)|=|\nabla \mathbb{G}\left[v_k( V(w_k+v)^p-w_k^p) \right](x)|
& \le& c_N(N-2)\int_{\R^N}\frac{h(y)}{|x-y|^{N-1}}dy
\\  &\le& c|x|^{-\theta_0-1}(1+|x|)^{2-N+\theta_0},
\end{eqnarray*}
thus, $|\mathcal{T} v |\in L^{q}(\R^N)$ for $\frac{N}{N-2}<q<\frac{N}{\theta_0};$ and $|\nabla \mathcal{T} v |\in L^{q}(\R^N)$ for $\frac{N}{N-1}<q<\frac{N}{\theta_0+1},$ where $\theta_0+1<N-1$.

For $\sigma\ge 1$, denote
$W^{1,\sigma}(\R^N)$ the Sobolev space with the norm
$$\norm{u}_{W^{1,\sigma}}=\left(\int_{\R^N} (|u|^\sigma+|\nabla u|^\sigma) dx\right)^{\frac1\sigma}.$$
Therefore, we see that $\mathcal{T} \mathcal{D}_\epsilon \subset W^{1,q_0}(\R^N)\cap \mathcal{D}_\epsilon$.\smallskip

{\it  We next show that the operator $\mathcal{T}$ is compact.} To this end, we only  have to prove that $W^{1,q_0}(\R^N)\cap \mathcal{D}_\epsilon$ is compact in $L^{q_0}(\R^N)$. Since the embedding $W^{1,q_0}(\R^N)\hookrightarrow L^{q_0}(\R^N)$ is locally compact in $\R^N$, letting $\{\zeta_j\}_j$ be a bounded functions in $W^{1,q_0}(\R^N)\cap \mathcal{D}_\epsilon$ with $\varepsilon>0$ and $\zeta\in L^p(\R^N)\cup \mathcal{D}_\epsilon$, then
for any $\eta>0$, there exist $R>0$, $j_\eta\in\N$ and a subsequence, still denote $\{\zeta_j\}_j$, such that for $j\ge j_\eta$,
$$\norm{\zeta_j-\zeta}_{L^p(B_R(0))}\le \frac\eta2\quad {\rm and}\quad \norm{\zeta_j}_{L^p(\R^N\setminus B_R(0))}+\norm{\zeta}_{L^p(\R^N\setminus B_R(0))}\le \frac\eta2,$$
therefore, we have that for $j\ge j_\eta$,
$$\norm{\zeta_j-\zeta}_{L^p(\R^N)}\le  \eta.$$
By the arbitrarily of $\eta$,  $W^{1,q_0}(\R^N)\cap \mathcal{D}_\epsilon\hookrightarrow L^{q_0}(\R^N)$ is compact and
we derive that $\mathcal{T}$ is a compact operator.

Observing that $\mathcal{D}_k$ is a closed and convex set in $L^{q_0}(\R^N)$, we now can apply Schauder
fixed point theorem to derive that there exists $v_k\in \mathcal{D}_\epsilon$ such that
$$\mathcal{T}v_k=v_k. $$
Since $|v_k(x)|\le  \epsilon |x|^{-\theta_0} (1+|x|)^{2-N+\theta_0}$, so $v_k$ is locally bounded in $\R^N\setminus\{0\}$,
then $v_k$ satisfies  (\ref{4.2}) by standard interior regularity results and $v_k$ is a  classical solution of (\ref{eq 3.1}).
 \hfill$\Box$
 \smallskip

 \begin{corollary}\label{Re 3.1}
 $(i)$  From the proof of Proposition \ref{pr 3.1}, the parameters $\epsilon$, $r$ could be fixed by
$$\epsilon=c_{\delta_0} k^p\quad{\rm and}\quad r=c_{\delta_0} k^{\frac{2p}{\tau_0}},$$
where $c_{\delta_0}>0$ is a constant depending on $\delta_0$.\smallskip

$(ii)$ Under the assumption of Proposition \ref{pr 3.1},  if $V\leq 1$, we can refine the solution  $v_k$ of (\ref{eq 3.1}) to be nonpositive,  derived in
 \begin{equation}\label{dom 1-}
\mathcal{D}_{\epsilon,-}:=\left\{w\in L^{q_0}(\R^N):\, - \epsilon|x|^{-\theta_0}(1+|x|)^{2-N+\theta_0}\leq w(x)\leq0 ,\ \forall\, x\in \R^N\setminus\{0\} \right\};
\end{equation}
if $  V\geq 1$, we can refine the solution  $v_k$ of (\ref{eq 3.1}) to be nonnegative,  derived in
 \begin{equation}\label{dom 1-}
\mathcal{D}_{\epsilon,+}:=\left\{w\in L^{q_0}(\R^N):\, 0\leq w(x)\leq  \epsilon|x|^{-\theta_0}(1+|x|)^{2-N+\theta_0},\ \forall\, x\in \R^N\setminus\{0\} \right\}.
\end{equation}
 \end{corollary}

\begin{remark}
In the critical case $p=p_c$, our construction of fast decaying solution
fails, due to the estimate of
$p\Gamma\ast (V w_k^{p-1}|v|)$ as the proof of Proposition \ref{pr 3.1},
where $\frac{p_cc_{p_c}^{p_c-1} }{\theta_0(N-2-\theta_0)}=1$ even with
$\theta_0$ taking the optimal value $\frac{N-2}2$ for $\theta_0(N-2-\theta_0)$, so there is no space for perturbing $V$ near origin.

\end{remark}

\subsection{Existence when $V$ is comparable to 1}

\begin{theorem}\label{teo 3.1}
Under assumptions of Theorem \ref{teo 1}, we let $V\geq1$.
Then there is $\nu_0\in(0,\, +\infty]$ such that for any $\nu\in(0,\,\nu_0)$, problem (\ref{eq 1.1}) has a $\nu$-fast decaying solution $u_\nu$,
which has singularity at the origin as (\ref{o1}) and the mapping $\nu\in (0,\,\nu_0)\mapsto  u_\nu$ is increasing, continuous and (\ref{1.3}) holds.

Moreover,
if (\ref{est 1}) holds for some $\alpha_1\geq0$ and $l_1>1$,
 then $\nu_0=+\infty$.

\end{theorem}
\noindent{\bf Proof.}   From Corollary \ref{Re 3.1}, we take $\epsilon=c_{\delta_0}k^p$
and let $k\in(0,\,k^*)$, then Proposition \ref{pr 3.1} implies that problem (\ref{eq 3.1})
has a nonnegative solution $v_{k}$ verifying (\ref{4.2}).
We denote
 $$\tilde u_{\nu_{k}}=w_{k}+v_{k} \geq w_k \quad{\rm and}\quad \tilde\nu_{k}=c_N\int_{\R^N}V\tilde u_{\nu_k}^pdx, $$
 then $\tilde u_{\nu_{k}}$ is a nonnegative classical solution of (\ref{eq 1.1}) such that
 \begin{equation}\label{5.2}
 \lim_{|x|\to0^+} \tilde u_{\nu_{k}}(x)|x|^{\frac{2}{p-1}}=c_p\quad{\rm and}\quad \lim_{|x|\to+\infty}\tilde u_{\nu_{k}}(x)|x|^{N-2}=\tilde \nu_{k},
 \end{equation}
where the second estimate is obtain by Lemma \ref{lm 1.3} and the fact that
$k\leq \tilde \nu_k\leq k+c_{\delta_0}k^p$.

To complete the proof, we divide into four steps.\smallskip

{\it Step 1. Existence by iteration method.}  We initiate from $v_0:= w_{k}$, denote by $v_n$ iteratively the unique solution of
 \begin{equation}\label{se 1}
  v_n=\Gamma\ast(Vv_{n-1}^p)\quad{\rm in}\quad\R^N\setminus\{0\},
  \end{equation}
that is,
$$
 \arraycolsep=1pt\left\{
\begin{array}{lll}
 \displaystyle -\Delta v_n=Vv_{n-1}^p\quad
    {\rm in}\quad   \R^N\setminus\{0\} ,\\[2mm]
 \phantom{    }
 \displaystyle  \lim_{|x|\to0} v_n(x)|x|^{N-2}=0.
 \end{array}\right.
$$
 As $-\Delta v_0= v_0^p$ in $\mathbb{R}^N\setminus\{0\}$ and $\displaystyle \lim_{|x|\to 0} v_0 |x|^{\frac{2}{p-1}}=c_{p}$,  by the Comparison Principle, we have that
$$v_1\ge v_0\quad{\rm in}\quad \R^N\setminus\{0\}.$$
Inductively, we can deduce that $v_n\ge v_{n-1}$ in $\R^N\setminus\{0\}$.
 Thus, the sequence $\{v_n\}_n$ is increasing.

Now we show that $\tilde u_{\nu_k}$ is an upper bound for $\{v_n\}_n$ for $k\in(0,\, k^*)$.
 We observe that $\tilde u_{\nu_k}$ is a solution of (\ref{eq 1.1}) and
$$Vw_k^p\leq V  \tilde u_{\nu_k}^p,\quad{\rm in}\ \ \R^N\setminus\{0\}.$$
Then Comparison Principle implies that
$$v_1\le  \tilde u_{\nu_k}  \quad{\rm in}\quad \R^N\setminus\{0\}.$$
Inductively, we see that for any $n\in\N$, we have that
$$v_n\le  \tilde u_{\nu_k} \quad{\rm in}\quad \R^N\setminus\{0\},$$
so $\{v_n\}_n$ has an  upper barrier  $\tilde u_{\nu_k}$.
Therefore, the sequence $\{v_n\}_n$ converges. Denote $\displaystyle  u_{\nu_k}:=\lim_{n\to\infty} v_n$, then for any compact set $K$ in $\R^N\setminus\{0\}$,
 and then $u_{\nu_k}$ verifies the equation $$-\Delta u= Vu^p\quad{\rm in}\quad K$$
and  then $u_{\nu_k}$ is a classical solution of (\ref{eq 1.1})   verifying
\begin{equation}\label{2.3+}
 w_k\le u_{\nu_k}\le \tilde u_{\nu_k} \quad{\rm in}\quad \R^N\setminus\{0\}.
\end{equation}
From the classification of isolated singularities of positive solutions to problem (\ref{eq 1.1}), we have that
$$\lim_{|x|\to0}u_{\nu_k}(x)|x|^{\frac{2}{p-1}}=c_{p}.$$
Here we let
 $$\nu_{k}=c_N\int_{\R^N}Vu_{\nu_k}^pdx.$$
 and then
 $$k\leq \nu_k\leq \tilde \nu_k\le k+c_{\delta_0}k^p\quad{\rm and}\quad
 \lim_{|x|\to+\infty}  u_{\nu_{k}}(x)|x|^{N-2}=  \nu_{k}$$
  hold by Lemma \ref{lm 1.3}.

 Thus,
\begin{equation}\label{3.5+}
\lim_{k\to0}\nu_k=0.
\end{equation}
Since $u_{\nu_k}\leq \tilde u_{\nu_k}\leq w_k+v_k$,
\begin{equation}\label{1.3--1}
 \lim_{k\to 0^+} \|w_k\|_{L^\infty_{loc}(\R^N\setminus\{0\})}=0,\quad   \lim_{k\to0^+} \|v_k\|_{L^\infty_{loc}({\R^N}\setminus\{0\})}=0,
 \end{equation}
so  $u_{\nu_k}$ verifies (\ref{1.3}). {\it Here and in what follows, we always denote $u_{\nu_k}$
the $\nu_k$-fast decaying solution of (\ref{eq 1.1}) derived by the sequence $v_n$ defined in (\ref{se 1}) with initial value $w_k$. }

\smallskip

{\it Step 2:  the mapping $k\mapsto \nu_k$ is increasing.} For $0<k_1<k_2$, by the increasing monotonicity of $w_k$, we have that $w_{k_1}<w_{k_2}$.
Let $\{v_{n,k_i}\}$ be sequence of (\ref{se 1}) with the initial data $v_0= w_{k_i}$, here $i=1,2$.

Let
$$\nu_{n,i}=\lim_{|x|\to+\infty} v_{n,k_i}(x)|x|^{N-2},\quad  i=1,2,\quad n=1,2,3,\cdots$$
We see that
$$\nu_{1,1}=c_N\int_{\R^N}Vw_{k_1}^p dx<c_N\int_{\R^N}Vw_{k_2}^p dx=\nu_{1,2} $$
and
\begin{eqnarray*}
\nu_{1,2}-\nu_{1,1}=c_N\int_{\R^N}V(w_{k_2}^p-w_{k_1}^p) dx\geq c_N\int_{\R^N}(w_{k_2}^p-w_{k_1}^p) dx
= k_2-k_1.
\end{eqnarray*}
Inductively, we have that for any $n\in\N$,
$$\nu_{n,2}- \nu_{n,1}\ge k_2-k_1,$$
which implies that the limit $u_{\nu_{k_1}}$ of $\{v_{n,k_1}\}$ and the limit $u_{\nu_{k_2}}$ of $\{v_{n,k_2}\}$ as $n\to+\infty$
verifies that
$$\lim_{|x|\to+\infty} u_{\nu_{k_2}}(x)|x|^{N-2}-\lim_{|x|\to+\infty} u_{\nu_{k_1}}(x)|x|^{N-2}\ge k_2-k_1,$$
 that is,
$$\nu_{k_2}-\nu_{k_1}\ge k_2-k_1.$$
As a conclusion, for any $k\in(0,\, k^*)$, there exists a $\nu_k>0$ such that
problem (\ref{eq 1.1}) has a solution $u_{\nu_k}$ such that
$$\lim_{|x|\to+\infty} u_{{\nu_k}}(x)|x|^{N-2}=\nu_k.$$

For $0<k_2\le k_1\le k_0$, then $w_{k_1,\mu}\ge w_{k_2,\mu}$  and
$v_{n,k_1}\ge v_{n,k_2}$ in $\R^N\setminus\{0\}$, so we have that $u_{\nu_{k_1}}\ge u_{\nu_{k_2}}$ in $\R^N\setminus\{0\}$.
That is to say  that the mapping $k\mapsto \nu_k$ is increasing.\smallskip

{\it Step 3. we prove that the mapping $k\in(0,\, k^*)\mapsto \nu_k$ is continuous.} Fix $\bar k\in (0,\, k^*)$ and $\delta<\frac12\min\{\bar k, k_0-\bar k\}$, then for $k\in(\bar k-\delta,\, \bar k+\delta)$ and $x\in\R^N\setminus\{0\}$, we have that
\begin{eqnarray*}
 &&| \bar w_{0}(-\ln|x|+b_{ 0}^{-1} \ln (\frac{k}{d_{0}})) - \bar w_{0}(-\ln|x|+b_{ 0}^{-1} \ln (\frac{\bar k}{d_{0}}))|
\\&& \le  b_{0}^{-1} |\ln k-\ln \bar k| \bar w_0'(-\ln|x|+b_{ 0}^{-1} \ln (\frac{\bar k+\delta}{d_{ 0}}))
\\&&\le c b_{ 0}^{-1} |\ln k-\ln \bar k| \bar w_0(-\ln|x|+b_{0}^{-1} \ln (\frac{\bar k+\delta}{d_{0}})),
\end{eqnarray*}
where the last inequality used (\ref{e 2.2}).

 For $|k-\bar k|<\delta$, we have that
 $$b_{0}^{-1} |\ln k-\ln \bar k| \le \bar c |k-\bar k|,  $$
where $\bar c= b_{ 0}^{-1}\ln2$.  So we have that
$$|v_{0,k}-v_{0,\bar k}|(x)\le  c|k-\bar k| w_{\bar k+\delta}(x)\leq c|k-\bar k| u_{\nu_{\bar k+\delta}}(x).  $$
Then by the increasing monotonicity of the mapping $k\mapsto \nu_k$, we have that
\begin{eqnarray*}
|v_{1,k}-v_{1,\bar k}|  \le  \mathbb{G}[V|v_{0,k}^p-v_{0,\bar k}^p|]
  &\le&  p\mathbb{G}[V v_{0,\bar k+\delta}^{p-1}|v_{0,k}-v_{0,\bar k}|] 
  \\&\le&  c |k-\bar k| \mathbb{G}[V u_{\nu_{\bar k+\delta}}^p]
=  c |k-\bar k| u_{\nu_{\bar k+\delta}}
\end{eqnarray*}
and since $v_{1,k},\, v_{1,\bar k}\leq u_{\nu_{\bar k+\delta}}$, it holds
\begin{eqnarray*}
|v_{2,k}-v_{2,\bar k}| \le \mathbb{G}[V|v_{1,k}^p-v_{1,\bar k}^p|]
\le  c  \mathbb{G}[Vu_{\nu_{\bar k+\delta}}^{p-1}|v_{1,k}-v_{1,\bar k}| ]
=  c |k-\bar k| u_{\nu_{\bar k+\delta}} ,
\end{eqnarray*}
inductively, we obtain that
$$
|v_{n,k}-v_{n,\bar k}|\leq c |k-\bar k| u_{\nu_{\bar k+\delta}},
$$
so, the following holds
$$|u_{\nu_{  k }}-u_{\nu_{\bar k}}|\le c_3|k-\bar k| u_{\nu_{\bar k+\delta}}\quad{\rm in}\ \  \R^N\setminus\{0\}.$$
Therefore, we have that
$$|\nu_k-\nu_{\bar k}|\le c\nu_{\bar k+\delta}|k-\bar k|. $$

Let $\displaystyle  \nu_0=\lim_{k\to k^*}\nu_{k}$,  then for any $\nu\in (0,\nu_0)$, problem (\ref{eq 1.1}) has a $\nu$-fast decaying solution $u_\nu$ verifying
\begin{equation}\label{origin}
 \lim_{|x|\to0}u_\nu(x)|x|^{\frac{2}{p-1}}=c_{p}.
\end{equation}

Observe that for any $\bar \nu\geq \nu_0$, if there is an upper barrier for the sequence (\ref{se 1}) with
the initial dada $w_k$, all above properties could be extended into $(0,\bar \nu)$.

Motivated by the fact that $\nu_k\geq k$ when $V\geq 1$,  let us denote
\begin{eqnarray*}
\nu_\infty &=& \sup\left\{\nu>0,\ {\rm there\ is}\ k\in(0,\,\nu]\ {\rm  problem\, (\ref{eq 1.1})\ has\ solution}\ u_{\nu}\geq w_k
\right.
\\&  &\left. \qquad \ {\rm in}\ \R^N\setminus\{0\} \   {\rm derived\ by\ the\ sequence\ (\ref{se 1})\ with\ the  \ initial\ data}\ w_k \right\}
\end{eqnarray*}
Note that using the above arguments, we can show that  the mapping $k\in(0,\, k_\infty)\mapsto \nu_k$ is increasing and continuous, where $k_\infty=\sup\{k>0: \nu_k<\nu_\infty\}$.\medskip

{\it Finally, we prove that $\nu_\infty=+\infty$ if (\ref{est 1}) holds for some $l_1>1$ and $\alpha_1\ge0$. }
By contradiction, we may assume that
 \begin{equation}\label{4.1+1}
 \nu_\infty<+\infty.
 \end{equation}
Now fix $\bar\nu\in(0,\nu_\infty)$ such that for $l_1>1$ and $\alpha_1\ge0$ in (\ref{est 1}) and $\bar \nu l_1^{N-2-\frac{2+\alpha_1}{p-1}} >\nu_\infty,$ where $N-2-\frac{2+\alpha_1}{p-1}>0$ by our assumption (\ref{v00}).
Denote
$$\psi_1(x)=l_1^{-\frac{2+\alpha_1}{p-1}} u_{\bar\nu}(l_1^{-1}x),\qquad\forall\, x\in\R^N\setminus\{0\}.$$
Let
$\bar k$ be the number such that $\nu_{\bar k}=\bar \nu$.
By direct computation, we have that
$$\lim_{|x|\to+\infty}\psi_1(x)|x|^{N-2}=\bar\nu  l_1^{N-2-\frac{2+\alpha_1}{p-1}}$$
and $-\Delta \psi_1=V_{l_1}\psi_1^p \quad {\rm in}\quad \R^N\setminus\{0\}, $
where $V_{l_1}(x)=l_1^{\alpha_1} V(l_1^{-1}x)\geq V(x)$ by (\ref{est 1}).

Note that $w_{l_1^{N-2-\frac{2+\alpha_1}{p-1}}\bar k}(x)=l_1^{-\frac{2+\alpha_1}{p-1}}w_{\bar k}(l_1^{-1}x)$ and
then we may initiate the iteration (\ref{se 1}) with  $v_0=w_{l_1^{N-2-\frac{2+\alpha_1}{p-1}}\bar k}$ and $\psi_{0}$
is an upper bound, so we have a solution $u_{\nu_{l_1^{N-2-\frac{2+\alpha_1}{p-1}}\bar k}}$ of (\ref{eq 1.1}) such that
$$  u_{\nu_{l_1^{N-2-\frac{2+\alpha_1}{p-1}}\bar k}}\ge  w_{l_1^{N-2-\frac{2+\alpha_1}{p-1}}\bar k}$$
that means
$$\nu_{l_1^{N-2-\frac{2+\alpha_1}{p-1}}\bar k}>\bar \nu l_1^{N-2-\frac{2+\alpha_1}{p-1}}>\nu_\infty, $$
which contradicts (\ref{4.1+1}).
So we have that $\nu_\infty=+\infty$.
The proof ends. \hfill$\Box$\medskip

\begin{theorem}\label{teo 3.2}
Under assumptions of Thorem \ref{teo 1}, we let $V\leq1$, i.e. $\beta=0$ and $c_\infty=1$ in the assumption $\mathcal{V}_0$ part $(ii)$.
Then there is   $\nu_0\in(0,+\infty]$ such that
for any $\nu\in(0,\,\nu_0)$, problem (\ref{eq 1.1}) has a $\nu$-fast decaying solution $u_\nu$,
which has singularity at the origin verifying (\ref{o1}) and the mapping $\nu\in (0,\,\nu_0)
 \mapsto  u_\nu$ is increasing, continuous and (\ref{1.3}) holds.

If (\ref{est 2}) holds for some $\alpha_2\leq0$ and $l_2>1$,
 then $\nu_0=+\infty$.
 \end{theorem}
\noindent{\bf Proof.} From Corollary \ref{Re 3.1}, we take $\epsilon=c_{\delta_0}k^p$,
and let $k\in(0,k^*)$, then Proposition \ref{pr 3.1} implies that problem (\ref{eq 3.1})
has a nonpositive solution $v_{k}$ verifying (\ref{4.2}).
Denote
 $$\tilde u_{\nu_{k}}=w_{k}+v_{k} \leq w_k \quad{\rm and}\quad \tilde\nu_{k}=\int_{\R^N}V\tilde u_{\nu_k}^pdx, $$
 then  $\tilde u_{\nu_{k}}$ is a nonnegative  classical solution of (\ref{eq 1.1}) such that
 \begin{equation}\label{5.2}
 \lim_{|x|\to0^+} \tilde u_{\nu_{k}}(x)|x|^{\frac{2}{p-1}}=c_p\quad{\rm and}\quad \lim_{|x|\to+\infty}\tilde u_{\nu_{k}}(x)|x|^{N-2}=\tilde \nu_{k},
 \end{equation}
where the second estimate is obtain by Lemma \ref{lm 1.3} and
$(k-c_{\delta_0}k^p)_+< \tilde \nu_k\leq k$.

Let $v_0:= w_{k}$ and denote $v_n$ iteratively
\begin{equation}\label{se 2}
v_n=\Gamma\ast(Vv_{n-1}^p)\quad{\rm in}\quad\R^N\setminus\{0\},
\end{equation}
 Since $V\leq 1$,  we have that the sequence $\{v_n\}_n$ is a decreasing sequence of functions. When $k\in(0, k^*)$, we can show that $\tilde u_{\nu_{k}}$, obtained in Proposition \ref{pr 3.1}, is a positive lower barrier for this sequence and we derive a  $\nu_k$-fast decaying solution $u_{\nu_k}$ of problem (\ref{eq 1.1}),
 which verifies (\ref{1.3}), where $\nu_{k}=c_N\int_{\R^N}Vu_{\nu_k}^pdx.$ When $k\geq k^*$, $\tilde u_{\nu_{k^*/2}}$ is alway a lower bound
 for the sequence $\{v_n\}_n$ and then a solution $u_{\nu_k}$ of problem (\ref{eq 1.1}).\smallskip

 {\it We show that the mapping $k\in(0,+\infty)\to \nu_k$ is increasing and continuous. }  For $0<k_1<k_2$, we have that $w_{k_1}<w_{k_2}$.
Let $\{v_{n,k_i}\}$ be sequence of (\ref{se 2}) with the initial data $v_{0,i}=w_{k_i}$, here $i=1,2$.

Let
$$\nu_{n,i}=\lim_{|x|\to+\infty} v_{n,k_i}(x)|x|^{N-2},\quad  i=1,2,\quad n=1,2,3,\cdots$$
We see that
$$\nu_{1,1}=c_N\int_{\R^N}Vw_{k_1}^p dx<c_N\int_{\R^N}Vw_{k_2}^p dx=\nu_{1,2} $$
and
$$0< \nu_{1,2}-\nu_{1,1}=c_N\int_{\R^N}V(w_{k_2}^p-w_{k_1}^p) dx\le c_N\int_{\R^N}(w_{k_2}^p-w_{k_1}^p) dx  = k_2-k_1.$$
Inductively,
we have that for any $n\in\N$,
$$0\le \nu_{n,2}- \nu_{n,1}\le k_2-k_1,$$
which implies that the limit $u_{\nu_{k_1}}$ of $\{v_{n,k_1}\}$ and the limit $u_{\nu_{k_2}}$ of $\{v_{n,k_2}\}$ as $n\to+\infty$
verifies that
$$0\le \lim_{|x|\to+\infty} u_{\nu_{k_2}}(x)|x|^{N-2}-\lim_{|x|\to+\infty} u_{\nu_{k_1}}(x)|x|^{N-2}\le k_2-k_1,$$
 that is,
$$0\le \nu_{k_2}-\nu_{k_1}\le k_2-k_1.$$
As a conclusion,  $k\to \nu_k$ is increasing and continuous.

Let
$$\nu_\infty=\lim_{k\to+\infty} \nu_k,$$
then we have that
for any $k\in(0,+\infty)$, there exists  $\nu_k\in (0,\nu_\infty)$ such that
problem (\ref{eq 1.1}) has a solution $u_{\nu_k}$ such that
$$\lim_{|x|\to+\infty} u_{{\nu_k}}(x)|x|^{N-2}=\nu_k.$$

{\it Finally, we prove that $\nu_\infty=+\infty$. } By contradiction, we may assume that
 \begin{equation}\label{4.1+}
 \nu_\infty<+\infty.
 \end{equation}
Now fix $\bar\nu\in(0,\nu_\infty)$, then there exist $\alpha_2\le 0$ and $l_2>1$ such that
 $$\bar \nu l_2^{N-2-\frac{2+\alpha_2}{p-1}} >\nu_\infty$$ and
 denote
$$\psi_2(x)=l_2^{-\frac{2+\alpha_2}{p-1}} u_{\bar\nu}(l_2^{-1}x),\quad\forall\, x\in\R^N\setminus\{0\}.$$
Let $\bar k$ be the number such that $\nu_{\bar k}=\bar \nu$.
By direct computation, we have that
$$\psi_2(x)\le l_2^{-\frac{2+\alpha_2}{p-1}}w_{\bar k}(l_2^{-1}x)),\quad\forall\, x\in\R^N\setminus\{0\} $$
and
$$-\Delta \psi_2=V_{l_2}\psi_2^p \quad {\rm in}\quad \R^N\setminus\{0\}, $$
where $V_{l_2}(x):=l_2^{\alpha_2}V(l_2^{-1}x)\leq V(x)$ by (\ref{est 2}).

Note that $w_{l_2^{N-2-\frac{2+\alpha_2}{p-1}}\bar k}(x)=l_2^{-\frac{2+\alpha_2}{p-1}}w_{\bar k}(l_2^{-1}x)$ and
then we may initiate the iteration (\ref{se 1}) with  $v_0=w_{l_2^{N-2-\frac{2+\alpha_2}{p-1}}\bar k}$ and $\psi_{0}$
is a lower bound, so we have a solution $u_{\nu_{l_2^{N-2-\frac{2+\alpha_2}{p-1}}\bar k}}$ of (\ref{eq 1.1}) such that
$$\psi_2\le u_{\nu_{l_2^{N-2-\frac{2+\alpha_2}{p-1}}\bar k}}\le  w_{l_2^{N-2-\frac{2+\alpha_2}{p-1}}\bar k}$$
that means
$$\nu_{l_2^{N-2-\frac{2+\alpha_2}{p-1}}\bar k}>\bar \nu l_2^{N-2-\frac{2+\alpha_2}{p-1}}>\nu_\infty, $$
which contradicts (\ref{4.1+}).
Thus, we have that $\nu_\infty=+\infty$.
  \hfill$\Box$\medskip

\subsection{ Proof of  main Theorems }

\noindent{\bf Proof of Theorem \ref{teo 1}.}
{\it Step 1: Existence and the decay at infinity.}  Let
$$V_1=1-(V-1)_-\quad{\rm and}\quad V_2=1+(V-1)_+,$$
where $a_\pm=\max\{0,\,\pm a\}$. Then $V_1,\,V_2$ are H\"{o}lder continuous and
$V=V_1V_2$ in $\R^N\setminus\{0\}$.

From Theorem \ref{teo 3.2},  there is $\nu_1\in(0,+\infty]$ such that for any $\nu\in(0,\,\nu_1)$, problem
\begin{equation}\label{eq 3.1-2}
 \left\{
\begin{array}{lll}
 \displaystyle  -\Delta  u= V_1 u^p\quad    &{\rm in}\quad   \R^N\setminus\{0\},\\[2mm]
 \displaystyle  \quad\ \ u>0\quad    &{\rm in}\quad   \R^N\setminus\{0\},
 \end{array}\right.
\end{equation}
  has a $\nu$-fast decaying solution $u_{\nu,1}$,
which has singularity at origin satisfying (\ref{o1})
 and the mapping $\nu\in (0,\,\nu_1)
 \mapsto  u_{\nu}$ is  increasing, continuous and (\ref{1.3}) holds.


We remark that for any $\nu\in (0,\, \nu_1)$, there exists a unique $k$ such that
$\nu=\nu_k$ and the solution $u_{\nu}$ is derived as the limit of the sequence
$$v_n=\Gamma\ast(V_1v_{n-1}^p)\quad{\rm in}\quad\R^N\setminus\{0\},$$
with the initial data $v_0= w_{k}$.  The mapping $k\mapsto \nu_k$ is increasing, continuous
and $\lim_{k\to0^+} \nu_k=0$.

As the proof of Theorem \ref{teo 3.1},  there is a $\mu_1\in(0,+\infty]$ such that    for any $\mu\in(0,\,\mu_1)$, problem
\begin{equation}\label{eq 3.1-1}
 \left\{
\begin{array}{lll}
 \displaystyle  -\Delta  u= V_2V_1 u^p\quad    &{\rm in}\quad   \R^N\setminus\{0\},\\[2mm]
 \displaystyle  \quad\ \ u>0\quad    &{\rm in}\quad   \R^N\setminus\{0\}
 \end{array}\right.
\end{equation}
  has a $\mu$-fast decaying solution $u_{\mu}$,
which has singularity (\ref{o1}) at the origin
 and the mapping $\mu\in (0,\,\mu_1)
 \mapsto  u_{\mu}$ is increasing, continuous and (\ref{1.3}) holds.


We remark that for any $\mu\in (0,\, \mu_1)$, there exists a unique $\nu\in(0,\,\nu_1)$ such that
$\mu=\mu_\nu$ and the solution $u_{\mu}$ is derived as the limit of the sequence
$$v_n=\Gamma\ast(V_2V_1v_{n-1}^p)\quad{\rm in}\quad\R^N\setminus\{0\}$$
with the initial data $v_0:= u_{\nu,1}$.  The map $\nu\mapsto \mu_\nu$ is increasing and continuous
and $\lim_{k\to0^+} \mu_\nu=0$.

As a consequence, for some $k^*\in(0,+\infty]$,  the map: $k\mapsto \nu_k\mapsto \mu_{\nu_k}\, $, denoting $\mu_k=\mu_{\nu_k}$, is continuous and increasing,    problem (\ref{eq 1.1})  has a  $\mu_k$-fast decaying  solution $u_{\mu_k}$,
which has singularity at the origin
$$
 \lim_{|x|\to0}u_\nu(x)|x|^{\frac{2}{p-1}}=c_{p},
$$
where $c_p$ is given in (\ref{cp}). From Proposition \ref{pr 3.1} and Corollary \ref{Re 3.1}, we have that
$$\mu_k\leq k+c_{\delta_0}k^p$$
which deduces  (\ref{1.3}).
 \hfill$\Box$\medskip

\noindent{\bf Proof of Theorem \ref{teo 2}. } Theorem \ref{teo 2} follows by Theorem \ref{teo 3.1} and Theorem \ref{teo 3.2} directly.\hfill$\Box$

   \setcounter{equation}{0}
\section{Existence of slow decay solution}

Under the assumption $(\mathcal{V}_1)$  part $(II)$, the mapping $\nu\in(0,\,\infty)\mapsto  u_\nu$ is increasing,
where $u_\nu$ is a $\nu$-fast decaying solution of problem (\ref{eq 1.1}), so our interest is to show the limit of $\{u_\nu\}_\nu$ as $\nu\to+\infty$ exists, denoting $\displaystyle u_\infty=\lim_{\nu\to+\infty} u_\nu$ if the limit exists, which is a very weak solution of (\ref{eq 1.1}) in the distributional sense that
$u_\infty\in L^1_{loc}(\R^N)\cap L^p_{loc}(\R^N,Vdx)$ satisfies the identity
\begin{equation}\label{id 1}
\int_{\R^N} u_\infty (-\Delta) \xi \,dx=\int_{\R^N} V u_\infty^p   \xi\,  dx,\quad\ \forall\, \xi\in C_c^\infty(\R^N).
\end{equation}

\begin{lemma}\label{lm 4.2}
 Let the hypotheses of Theorem \ref{teo 3} hold,
then $u_\nu$ is radially symmetric and decreasing with respective to $|x|$.

\end{lemma}
{\bf Proof.} Since $V$ is radially symmetric and decreasing with respect to $|x|$ and
$u_\nu$ decays as $\nu|x|^{2-N}$ at infinity, so it is available to apply the moving planes method, see \cite{BN,GNN}
to obtain that $u_\nu$ is the radially symmetric and decreasing  with respect to $|x|$.\hfill$\Box$

\begin{lemma}\label{lm 4.1+}
Let the hypotheses of Theorem \ref{teo 3} hold, then there exists $c>0$ independent of $\nu$ such that
$$ \norm{u_\nu}_{L^1_{loc}(\R^N)}\le c \quad{\rm and}\quad \norm{u_\nu}_{L^p_{loc}(\R^N,Vdx)}\le c. $$

\end{lemma}
{\bf Proof.} Let
$$U_0(x)=\frac{c_N}{(1+|x|^2)^{\frac{N-2}{2}}},$$
then
$$-\Delta U_0=U_0^{2^*-1}\quad {\rm in}\quad \R^N,$$
where $2^*=\frac{2N}{N-2}$, $c_N=(N(N-2))^{\frac{N-2}{4}}$.

For $\epsilon\in(0,\frac14)$, denote
 \begin{equation}\label{eta}
 \eta_\epsilon(x)=\eta_0(\epsilon|x|),\quad x\in\R^N,
\end{equation}
where $\eta_0:[0,+\infty)\to[0,1]$ is a smooth increasing function such that
$$\eta_0(t)=0,\quad \forall\,t\ge 2\quad {\rm and}\quad \eta_0(t)=1,\quad  \forall\,t\in[0,1].$$

Take $U_0\eta_\epsilon^2$ as a test function of (\ref{id 1}), then by H\"{o}lder inequality, we have that
\begin{eqnarray*}
\int_{\R^N} V u_\nu^p   U_0  \eta_\epsilon^2 dx &=&\int_{\R^N} u_\nu (-\Delta) (U_0\eta_\epsilon^2) dx
\\&=&\int_{\R^N} u_\nu   (U_0^{2^*-1}\eta_\epsilon^2+4\eta_\epsilon \nabla U_0\cdot\nabla \eta_\epsilon +U_0(-\Delta)(\eta_\epsilon^2) )dx
\\&\le&  \left(\int_{\R^N} V u_\nu^p   U_0\eta_\epsilon^2  dx\right)^{\frac1p}  \left(\int_{B_{\frac2\epsilon}} V^{1-p}  U_0^{1-p+(2^*-1)\frac{p}{p-1}}   dx\right.
\\&&+ \epsilon \norm{\eta_0}_{C^1(\R)}\int_{B_{\frac2\epsilon}(0)\setminus B_{\frac1\epsilon}(0)} V^{1-p}  U_0^{1-p+(2^*-1)\frac{p}{p-1}} |\nabla U_0|^{p-1}  dx
\\&&\left.+
\epsilon^2 \norm{\eta_0}_{C^2(\R)}\int_{B_{\frac2\epsilon}(0)\setminus B_{\frac1\epsilon}(0)} V^{1-p}  U_0^{p-1}   dx\right)^{1-\frac1p}
\\&=& c_0 \left(\int_{\R^N} V u_\nu^p   U_0 \eta_\epsilon^2  dx\right)^{\frac1p},
\end{eqnarray*}
where $c_0>0$ is dependent on $\epsilon,V$ but it is independent of $\nu$.

So we have that
\begin{eqnarray*}
 c_N(1+\frac1{\epsilon^2})^{-\frac{N-2}2}\int_{B_{\frac1\epsilon}(0)} V u_\nu^p    dx \le   \int_{B_{\frac1\epsilon}(0)} V u_\nu^p   U_0  dx
    \le\int_{\R^N} V u_\nu^p   U_0  \eta_\epsilon^2 dx \le c_0^{\frac{p}{p-1}},
\end{eqnarray*}
that is,
$$\norm{u_\nu}_{L^p_{loc}(\R^N,Vdx)}\le c $$
for some $c>0$ independent of $\nu$.

Furthermore,
$$\int_{\R^N} u_\nu   U_0^{2^*-1}\eta_\epsilon^2 dx \le c_0 \left(\int_{\R^N} V u_\nu^p   U_0 \eta_\epsilon^2 dx\right)^{\frac1p}\le c_0^{\frac{p}{p-1}}$$
and
$$\norm{u_\nu}_{L^1_{loc}(\R^N)}\le c. $$
The proof ends.\hfill$\Box$\medskip

\noindent{\bf Proof of Theorem \ref{teo 3}.}
From Theorem \ref{teo 2} and Lemma \ref{lm 4.1+}, the mapping $\nu\in(0,\,\infty)\mapsto  u_\nu$ is increasing and uniformly bounded in
$L^1_{loc}(\R^N)\cap L^p_{loc}(\R^N,Vdx)$, so there exists $u_\infty\in L^1_{loc}(\R^N)\cap L^p_{loc}(\R^N,Vdx)$ such that
$$u_\nu\to u_\infty \quad {\rm as}\quad \nu\to+\infty \quad{\rm a.e.\ in }\ \Omega \ \ {\rm and\ in }\  L^1_{loc}(\R^N)\cap L^p_{loc}(\R^N,Vdx). $$

It is known in \cite{L} that $u_\nu$ is also a weak solution of (\ref{eq 1.1}), i.e.
\begin{equation}\label{id 2}
\int_{\R^N} u_\nu (-\Delta) \xi dx=\int_{\R^N} V u_\nu^p   \xi  dx,\qquad\forall\, \xi\in C_c^\infty(\R^N).
\end{equation}
Passing to the limit of (\ref{id 2}), we obtain that $u_\infty$ is a weak solution of (\ref{eq 1.1})
in the sense of (\ref{id 1}).

Furthermore, by Lemma \ref{lm 4.2}  we have that $u_\nu$ is radially symmetric and decreasing  with respect to $|x|$, so is $u_\infty$.
So we have that  $u_\infty\in L^\infty_{loc}(\R^N\setminus \{0\})$, then $Vu_\infty^p$ is in $L^\infty_{loc}(\R^N\setminus \{0\})$. By standard
regularity results, we have that $u$ is a classical solution of (\ref{eq 1.1}).

Since $u_\nu$ verifies (\ref{o1}) at the origin for any $\nu>0$ and $u_\infty$ is the limit of an increasing sequence $\{u_\nu\}_\nu$,
then we have that
\begin{equation}\label{4.o1}
 \liminf_{|x|\to0}u_\infty(x)|x|^{\frac{2}{p-1}}\ge c_{p}.
\end{equation}
Next we claim
\begin{equation}\label{4.o2}
 \limsup_{|x|\to0}u_\infty(x)|x|^{\frac{2}{p-1}}\le c_{p},
\end{equation}
From \cite[Theorem 2.1]{NS}, there exists $c>0$ such that
$$u_\infty(x)\le c|x|^{-\frac2{p-1}},\quad x\in B_1(0)\setminus\{0\}.$$
And by \cite[Theorem B]{A}, we have that (\ref{4.o2}) and (\ref{o1}) holds true for $u_\infty$.

Finally, we claim (\ref{1.3.2}).  Denote by
$$  v(x)=u_\infty(\frac{x}{|x|^2})\quad{\rm for}\quad x\in\R^N\setminus\{0\}.$$
By direct computation, we have that
$$\nabla v(x)=\nabla u_\infty\left(\frac{x}{|x|^2}\right)\frac1{|x|^2}-2\left(\nabla u_\infty(\frac{x}{|x|^2})\cdot x\right)\frac{x}{|x|^4}$$
and
$$\Delta v(x)=\frac1{|x|^4} \Delta u_\infty(\frac{x}{|x|^2})+\frac{2(2-N)}{|x|^4} \left(\nabla u_\infty(\frac{x}{|x|^2})\cdot x\right).$$
Let $u^\sharp(x)=|x|^{2-N}v(x)$,  then for $x\in \R^N\setminus\{0\}$, we get that
\begin{eqnarray*}
 -\Delta u^\sharp(x) = -\Delta v(x) \,|x|^{2-N} -2\nabla v(x)\cdot(\nabla |x|^{2-N})
 = |x|^{-2-N} (-\Delta )u(\frac{x}{|x|^2})
  = V^\sharp(x)u^\sharp(x)^p ,
\end{eqnarray*}
where
$$ V^\sharp(x)= |x|^{-2-N+p(N-2)}V(\frac{x}{|x|^2}). $$

We see that
$$ -\Delta u^\sharp(x)= V^\sharp(x)u^\sharp(x)^p\quad {\rm in}\quad \R^N\setminus\{0\},$$
where
$$V^\sharp(x)\sim |x|^{\varrho}\quad{as}\quad |x|\to0\quad\quad{\rm with}~~\varrho=-\alpha-4+(p-1)(N-2).$$
Note that it implies by  (\ref{poten 1}) that
 $\varrho\in(-2,0)$.

\medskip

Now claim that $p\in\Big(\frac{N+\varrho}{N-2},\, \frac{N+2+\varrho}{N-2}\Big)$. Note that
$p>\frac{N+\varrho}{N-2}$ follows by the fact $p>\frac{N}{N-2}$ and
$p<\frac{N+2+\varrho}{N-2}$ is equivalent to
$(N-2)p-N-2< -4+(p-1)(N-2)-\alpha$,
which is true thanks to $\alpha\leq 0$.

Then by  \cite[Theorem 2.1]{NS} and \cite[Theorem 3.3]{GSG} we have that
$$\frac1c|x|^{-\frac{2+\varrho}{p-1}}\le u^\sharp(x)\le c|x|^{-\frac{2+\varrho}{p-1}},\quad\forall\, x\in B_1(0)\setminus\{0\},$$
where $c>1$ and
$$-\frac{2+\varrho}{p-1}=-(N-2)+\frac{2+\alpha}{p-1}.$$
By Kelvin transformation, we turn back that
$$\frac1c|x|^{-\frac{2+\alpha}{p-1}}\le u_\infty(x)\le c|x|^{-\frac{2+\alpha}{p-1}},\quad x\in\R^N\setminus B_1(0)$$
for some $c>1$.\hfill$\Box$

\medskip

\bigskip

  \noindent{\bf Acknowledgements:} H. Chen  is supported by  Natural Science Foundation  of China (No:11726614, 11661045).
 X. Huang is supported by NSF of China, No.11701181.
F. Zhou is supported by the Natural Science Foundation  of China (No:11726613, 11431005),
by Science and Technology Commission of Shanghai Municipality (No:18dz2271000).


\begin{thebibliography}{99}



\bibitem {AS} S.N. Armstrong and B. Sirakov, Sharp Liouville results for fully
nonlinear equations with power-growth nonlinearities,  Ann.
della Scuola Normale Super. di Pisa Classe di scienze 10 (2011),
711-728.


\bibitem{A} P. Aviles, Local behaviour of the solutions of some elliptic equations,
  Comm. Math. Phys. 108 (1987), 177-192.


\bibitem{BN}
H. Berestycki and L. Nirenberg, On the method of moving planes and
the sliding method,  Bol. Soc. Brasileira Mat. 22  (1991), 1-37.


\bibitem{B} M. Bidaut-V\'{e}ron, Local behaviour of solutions of a class of nonlinear elliptic systems,
  Adv. Diff.  Eq. 5  (2000),  147-192.

 \bibitem{BP} M. Bidaut-V\'{e}ron and S. Pohozaev, Nonexistence results and estimates
for some nonlinear elliptic problems,  J. Anal. Math. 84 (2001), 1-49.




\bibitem{BV}  H. Brezis  and L. V\'{e}ron,
 Removable singularities for some nonlinear elliptic equations,
  Arch.  Ration. Mech. Anal. 75 (1980), 1-6.

 

 \bibitem {CGS} L.  Caffarelli, B. Gidas and J. Spruck, Asymptotic symmetry and local behaviour of semilinear
elliptic equation with critical Sobolev growth, Comm. Pure Appl. Math. 42 (1989), 271-297.


 \bibitem {CL} C. Chen and  C. Lin, Existence of positive weak solutions with a prescribed singular set of
semi-linear elliptic equations,   J. Geom. Anal. 9(2) (1999), 221-246.



\bibitem {C} H. Chen, S. Alhomedan, H. Hajaiej, P. Markowich,
Fundamental solutions for Schr\"{o}dinger operators with general inverse square potentials,
  Appl. Anal. 97 (2018), 787-810.

  \bibitem {CFY}  H. Chen, P. Felmer and J. Yang, Weak solutions of semilinear elliptic equation involving Dirac mass,
   Ann. I. H.  Poincar\'{e}-AN 35 (2018),  729-750.

 \bibitem {CQZ} H. Chen A. Quaas and F. Zhou,
On nonhomogeneous elliptic equations with the Hardy-Leray potentials, Accepted by   J. Anal.  Math.,  arXiv:1705.08047. 



 \bibitem {CPZ} H. Chen, R. Peng and F. Zhou, Nonexistence of positive supersolution to a class of semilinear elliptic equations and systems in an exterior domain,  Sci. China Math.  (2019), Doi.org/10.1007/s11425-018-9447-y.

\bibitem{DDG} E.  Dancer, Y. Du and Z. Guo,   Finite Morse index solutions of an elliptic equation with supercritical exponent,
  J. Diff. Eq. 250 (2011), 3281-3310.

 \bibitem{DPM} J. Davila, M. Del Pino,  M. Musso and  J. Wei,
  Fast and slow decay solutions for supercritical elliptic problems in exterior domains,
   Calc. Var.   Part. Diff. Eq. 32 (2008), 453-480.

 \bibitem{DLY} Y. Deng, Y. Li and F. Yang,
  A note on the positive solutions of an inhomogeneous elliptic equation on $\R^N$,
  J.  Diff. Eq. 246 (2009), 670-680.

\bibitem{GNN} B. Gidas, W-M. Ni and L. Nirenberg, Symmetry and related properties via the maximum principle,   Comm. Math. Phys. 68 (1979), 209-243.  
 

 \bibitem{GSG} B. Gidas and J. Spruck,   Global and local behavior of positive solutions of nonlinear elliptic
equations,   Comm. Pure Appl. Math. 42  (1989)  271-297.



\bibitem{L2}   Y. Li,  On the positive solutions of the Matukuma equation,
  Duke Math. J. 70 (1993),  575-589.


\bibitem {L}  P. Lions,   Isolated singularities in semilinear problems,  J. Diff. Eq. 38 (1980),   441-450.



\bibitem {MP}   R. Mazzeo and F. Pacard, A construction of singular solutions for a semilinear elliptic equation using asymptotic analysis,
 J. Diff. Geom. 44 (1996),  331-370.

\bibitem {MP1}  R. Mazzeo and F. Pacard, Constant scalar curvature metrices with isolated singularities,
 Duke Math. J. 99 (1999), 353-418.

\bibitem {N}  W-M. Ni, Some aspects of semilinear elliptic equations in $R^N$, in "Nonlinear Diffusion Equations and Their Equilibrium States", Springer-Verlag, New York, 1988.

\bibitem {NS} W-M. Ni and J. Serrin, Nonexistence theorems for singular solutions of quasilinear partial qquations, Comm. Pure Appl. Math. 39 (1986), 379-399.










\bibitem {GV} A. Gmira  and L. V\'{e}ron, Boundary singularities of solutions of some
nonlinear elliptic equations, Duke Math. J.  64 (1991), 271-324
.

\bibitem {VV} J. Vazquez  and L. V\'{e}ron,
 Removable singularities of some strongly nonlinear elliptic equations,
  Manuscripta Mathematica 33 (1980), 129-144.

 
\bibitem {MV2} M. Marcus  and L. V\'{e}ron, The boundary trace of positive
solutions of semilinear elliptic equations: the supercritical case,
  J. Math. Pures Appl. 77 (1998), 481-524.

\bibitem {MV3} M. Marcus  and L. V\'{e}ron, Removable singularities and boundary
traces,  J. Math. Pures Appl. 80 (2001), 879-900.

\bibitem {MV4} M. Marcus  and L. V\'{e}ron, The boundary trace and generalized B.V.P. for
semilinear elliptic equations with coercive absorption,   Comm.
Pure Appl. Math. 56 (2003), 689-731.

 \bibitem {V0}  L. V\'{e}ron, Singular solutions of some nonlinear elliptic equations,   Nonlinear Anal. T. M. $\&$ A. 5 (1981),  225-242.

\bibitem {V}  L. V\'{e}ron, Singularities of solutions of second order quasilinear equations,
  Pitman Research Notes in Mathematics Series Vol 353,    1996.


\end{thebibliography}
\end{document}